\documentclass[11pt]{article}
\usepackage{epsfig,amssymb}

\begin{document}

\title{
On implicit ODEs with hexagonal web of solutions}

\author{{\Large Agafonov S.I. } \\
\\
Institut f\"{u}r Mathematik\\
Martin-Luther-Universit\"{a}t
Halle-Wittenberg \\ D-06099 Halle (Saale), Germany\\
e-mail: {\tt sergey.agafonov@mathematik.uni-halle.de} }
\date{}
\maketitle

\unitlength=1mm

\newtheorem{theorem}{Theorem}
\newtheorem{proposition}{Proposition}
\newtheorem{lemma}{Lemma}
\newtheorem{corollary}{Corollary}
\newtheorem{definition}{Definition}
\newtheorem{example}{Example}

\pagestyle{plain}

\begin{abstract}
\noindent Solutions of an implicit ODE form a web. Already for
cubic ODEs the 3-web of solutions has a nontrivial local
invariant, namely the curvature form. Thus any local
classification of implicit ODEs necessarily has functional moduli
if no restriction on the class of ODEs is imposed. Here the most
symmetric case of hexagonal 3-web of solutions is discussed, i.e.
the curvature is supposed to vanish identically. A finite list of
normal forms is established under some natural regularity
assumptions. Geometrical meaning of these assumptions is that the
surface, defined by ODE in the space of 1-jets,
 is smooth as well as the criminant, which is the critical set of
this surface projection to the plane.
(secondary).
\end{abstract}

\section{Introduction}
Consider an implicit ordinary differential equation
\begin{equation}\label{implicit} F(x,y,p)=0
\end{equation} with a smooth or real analytic $F$. This ODE  defines a surface $M$:
\begin{equation}\label{surface}
M:=\{(x,y,p)\in {\mathbb R ^2 \times \mathbb P}^1(\mathbb R): \
F(x,y,p)=0\},
\end{equation}
where $(x,y,p)$ are coordinates in the jet space $J^1(\mathbb R,
\mathbb R)$ with $p=\frac{dy}{dx}$. Generically the condition
$\left. {\rm grad} (F)\right |_{F(x,y,p)=0}\ne 0$
 holds true for any point $m=(x,y,p)\in M$, i.e.  $M$ is smooth.
If the projection $\pi : M\to \mathbb R^2,\  (x,y,p) \mapsto
(x,y)$ is a local diffeomorphism at a point $m\in M$ then this
point is called {\it regular}. In some neighborhood of the
projection of a regular point $\pi (m)$  equation (\ref{implicit})
can be solved for $p$ thus defining an explicit ODE.

If the projection $\pi $ is not a local diffeomorphism at $m$,
then the point $m$ is called a {\it singular} point of implicit
ODE (\ref{implicit}). The set of all singular points is called the
{\it criminant} of equation (\ref{implicit}) or the {\it apparent
contour} of the surface $M$ and will be denoted by $C$:
\begin{equation}\label{criminant}
C:=\{(x,y,p)\in \mathbb R^2 \times  \mathbb P^1(\mathbb R):
F(x,y,p)= F_p(x,y,z)=0 \},
\end{equation}
where the low subscript denotes a partial derivative:
$F_p=\frac{\partial F}{\partial p}$.

Studying of generic singular points of implicit ODEs was initiated
by Thom in \cite{Te}. Due to Whitney's Theorem such points are
folds and cusps of the projection $\pi$.  Local normal forms for
generic singularities were conjectured by Dara in \cite{Ds} and
for a generic fold point were established by Davydov in \cite{Dn}.
The classification list for a generic fold point of the projection
$\pi$ is exhausted by a well folded saddle point, a well folded
node point, a well folded focus point and the {\it regular}
singular point, where contact plane is transverse to the
criminant. Cusp points were studied by Dara \cite{Ds}, Bruce
\cite{B} and Hayakawa, Ishikawa, Izumiya, Yamaguchi \cite{HIIY}.
\\Usually the following {\it regularity condition} is imposed at
each point of the criminant:
\begin{equation}\label{regularity}
{\rm rank}((x,y,p)\mapsto (F,F_p))=2.
\end{equation}
 This regularity condition
implies that the criminant is a smooth curve. At each point $m$
outside the criminant $C$ the contact plane $dy-pdx=0$ cuts the
tangent plane $T_mM$ along a line thus giving a direction field,
which takes the form:
\begin{equation}\label{direction_field}
\tau=[F_p:pF_p:-(F_x+pF_y)],
\end{equation}
in the coordinates $(x,y,p)$. This direction field is called the
{\it characteristic field} of $M$. The projection $\pi (\gamma )$
of an integral curve $\gamma \subset M$ of the characteristic
field $\tau$ is called a {\it solution} of ODE (\ref{implicit}).
If $F_{pp}\ne 0$ at a point $m\in C$ on the criminant then
(\ref{implicit}) reduces locally to quadratic in $p$ (or {\it
binary}) ODE. Such equations were the subject of intensive study.
See, for example, \cite{BTb},\cite{BTd},\cite{Ci},\cite{GOTp}.

 Suppose equation (\ref{implicit}) has a triple
root $p_0$ at $(x_0,y_0)$  then the equation $F=0$ can be written
locally as a cubic equation
\begin{equation}\label{general_cub}
p^3+a(x,y)p^2+b(x,y)p+c(x,y)=0,
\end{equation}
as follows from the Division Theorem. Thus if in a domain
$U\subset \mathbb R^2$ outside the discriminant curve $\Delta:
=\pi (C)$ this cubic equation has 3 real roots $p_1,p_2,p_3$, we
have 3-web formed by solutions of (\ref{implicit}). A generic
3-web has a nontrivial invariant. In differential-geometric conext
this invariant is the curvature form of the web. Therefore any
general local classification of cubic implicit ODEs
(\ref{general_cub}) necessarily has functional moduli (cf.
\cite{Dn}). Moreover, this invariant is topological in nature
hence even topological classification will have functional moduli
if no restriction is imposed on the class of ODE. (See also
\cite{Nv} and \cite{Ne}, where web structure was used for studying
geometric properties of differential equations.)

In this paper we consider cubic ODEs (\ref{general_cub}) with a
hexagonal web of solutions. Equations of this type describe, for
example, webs of characteristics on solutions of {\it integrable}
systems of three PDEs of hydrodynamic type (see \cite{Fc},
\cite{Fi}). Another example is WDVV associativity equation (see
Example \ref{associativ} below).

\begin{definition}\label{def_hex}
Let $U\subset \mathbb R^2$ be the open set, where
(\ref{general_cub}) has 3 real roots $p_1,p_2,p_3$ and suppose
$U\ne \emptyset$. We say that (\ref{general_cub}) has a hexagonal
3-web of solutions if for the projection $\pi (m)$ of each regular
point $m\in M$ with $\pi (m)\in U$ there is a local diffeomorphism
at $\pi(m)$ mapping the solutions of (\ref{general_cub}) to three
families of parallel lines.
\end{definition}
The case of hexagonal web of solutions is also the most symmetric,
i.e. the Lie symmetry pseudogroup of (\ref{implicit}) at a regular
point has the largest possible dimension 3. The list of normal
forms turnes out to be finite provided regularity condition
(\ref{regularity}) is satisfied. These forms are given by the
following examples.

\begin{figure}[th]
\begin{center}
\epsfig{file=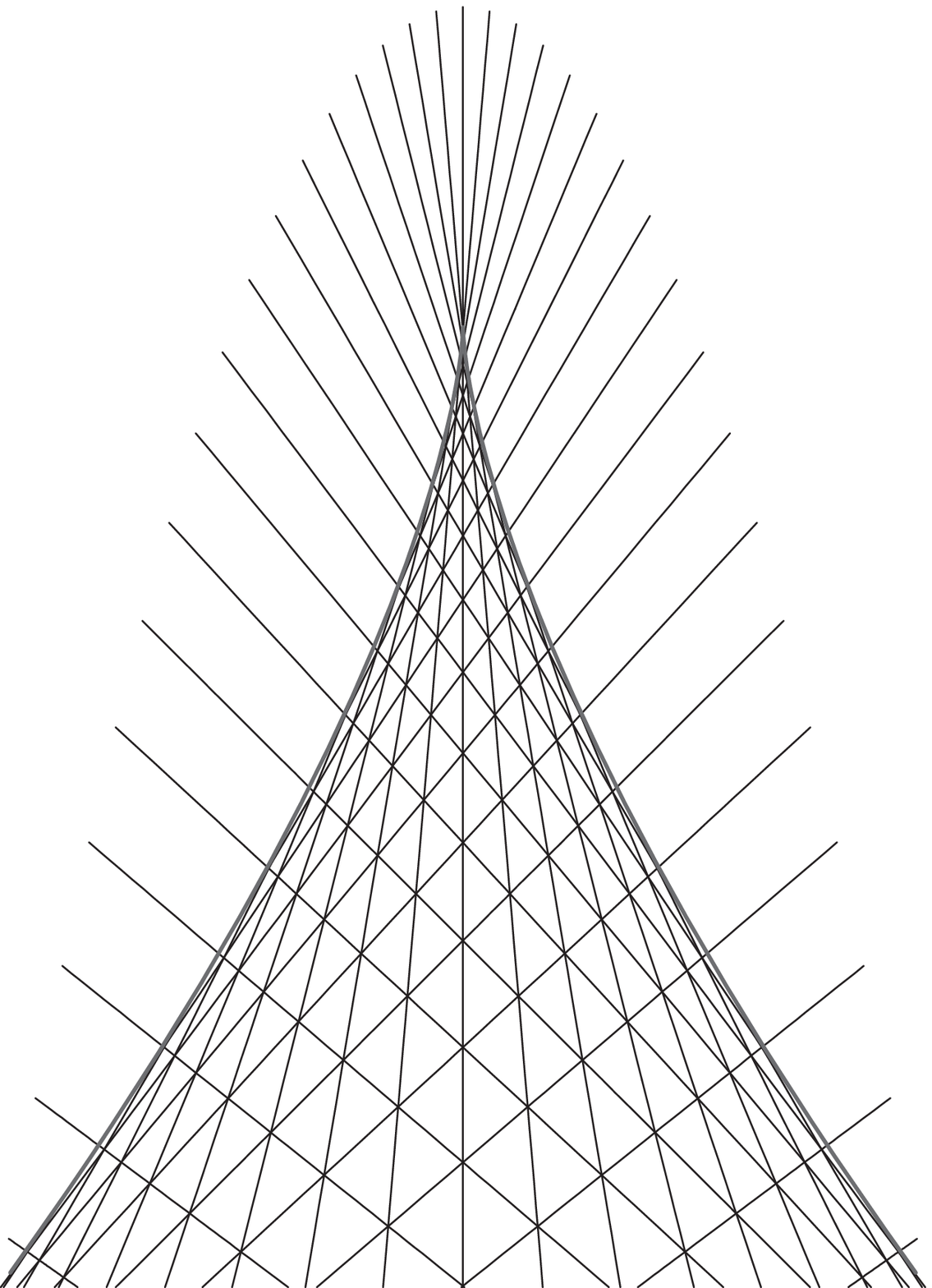,width=40mm} \ \ \ \ \ \ \ \ \ \ \ \ \ \ \ \
\ \ \ \ \ \ \
\epsfig{file=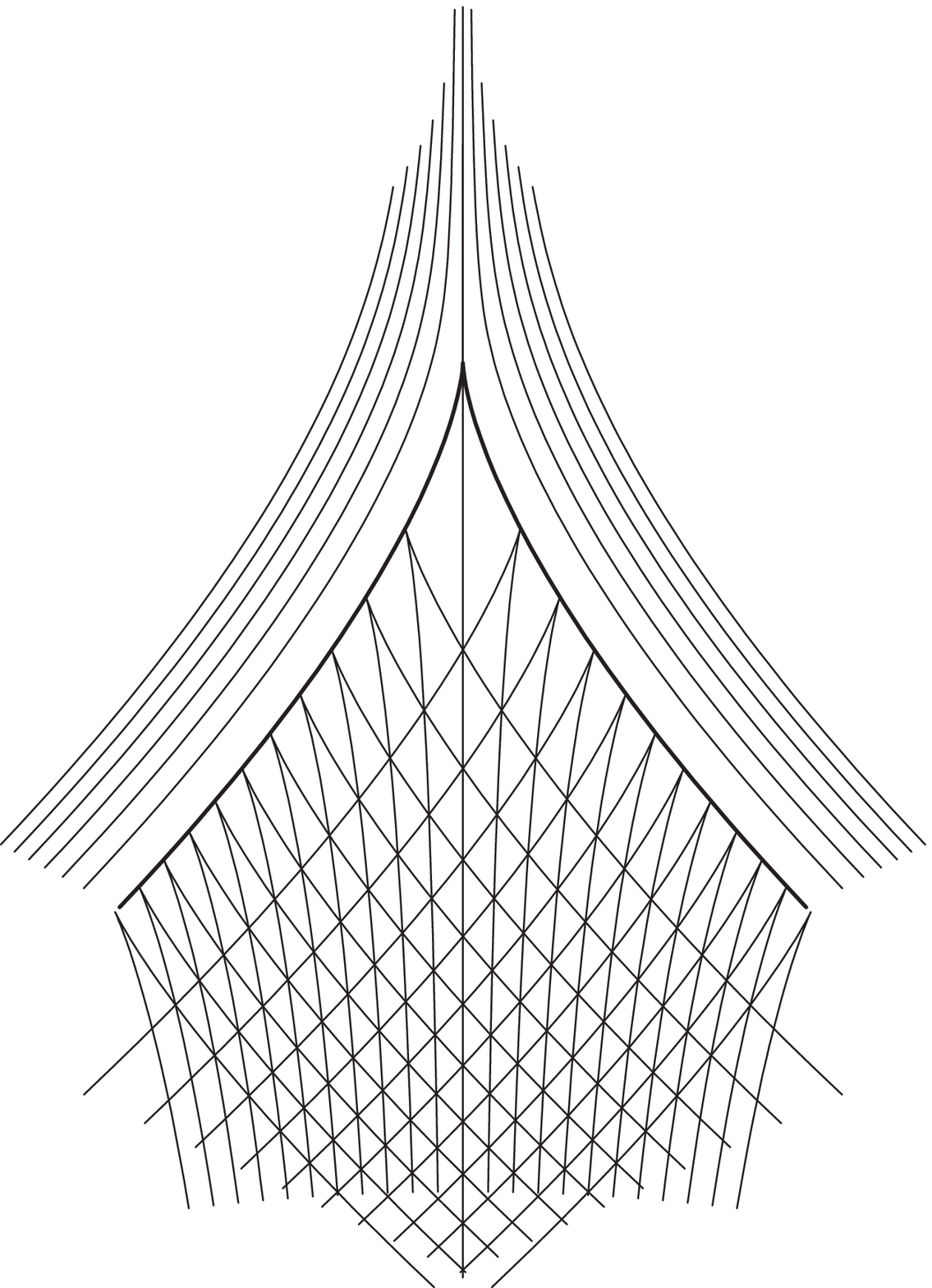,width=40mm} \caption{Solutions
of $p^3+px-y=0$ (left) and $p^3+2xp+y=0$ (right) with horizontal
y-axis.} \label{Pic_sol}
\end{center}
\end{figure}

\begin{example} {\rm The classical Graf and Sauer theorem
\cite{GSg} claims that a 3-web of straight lines is hexagonal iff
the web lines are tangents to an algebraic curve of class 3, i.e.
the dual curve is cubic. This implies immediately that the
following cubic Clairaut equation has a hexagonal 3-web of
solutions:
\begin{equation}\label{Clairaut_0}
p^3+px-y=0.
\end{equation}
The solutions are the lines $p=const$ enveloping a semicubic
parabola. (See Fig. \ref{Pic_sol}) Note that the contact plane is
tangent to $M$ along the criminant, i.e. the criminant is a {\it
Legendrian} curve.}
\end{example}

\begin{example}\label{associativ}  {\rm  Consider  associativity equation
$$
u_{xxx}=u^2_{xyy}-u_{xxy}u_{yyy},
$$
describing 3-dimensional Frobenius manifolds (see  \cite{Df}).
Each of its solutions $u(x,y)$ defines a characteristic web in the
plane, which is hexagonal as was shown by Ferapontov \cite{Fh}.
Characteristics are integral curves of the vector field
$$
\partial _x-\lambda (x,y)\partial_y,
$$
where $\lambda$  satisfy the  characteristic equation
$$
\lambda ^3+u_{yyy}\lambda ^2-2u_{xyy}\lambda + u_{xxy}=0.
$$
For the solution $u=\frac{x^2y^2}{4}+\frac{x^5}{60}$  the
characteristic equation becomes
\begin{equation}\label{Clairaut}
p^3+2xp+y=0
\end{equation}
after the substitution $x \to -x$, $y\to -y$, $\lambda\to -p$. The
criminant of this ODE is not Legendrian and the solutions have
ordinary cusps on the discriminant (see Fig. \ref{Pic_sol}). The
discriminant  is also a solution.  In the analytic setting the
above two normal forms were conjectured by Nakai in
\cite{Nw}.}\end{example}
We find also local normal forms at points, where the projection
$\pi$ has a fold, i.e. the cubic ODE factors out to a quadratic
and a linear terms.

\begin{example} {\rm Suppose the criminant of ODE is Legendrian then this ODE
 is locally equivalent to
\begin{equation}\label{q1}
p^2=y.
\end{equation}
Solutions of this ODE together with the lines $dx=0$ form a
hexagonal 3-web (see Fig. \ref{Pic_sol}). In fact, both the lines
$dy=0$ or the parabolas $2dy-xdx=0$ also supplement the 2-web of
solutions of (\ref{q1}) to a hexagonal 3-web, but the surfaces of
the corresponding cubic equations $p(p^2-y)=0$ and
$(2p-x)(p^2-y)=0$ are not smooth at $m=(0,0,0)$. If we agree to
consider a quadratic equation as a cubic with one root at
infinity, than equation (\ref{q1}) is the third normal form in our
list. The following coordinate change
$y=\tilde{y}+\frac{\tilde{x}^2}{4}$, $x=\tilde{x}$ straightens the
solutions, transforming ODE (\ref{q1}) to a quadratic Clairaut
equation
$$
\tilde{p}^2+\tilde{p}\tilde{x}-\tilde{y}=0.
$$
As the the lines $dx=0$ are preserved this example is also a
special case of the Graf and Sauer Theorem.

 }\end{example}

\begin{figure}[th]
\begin{center}
\epsfig{file=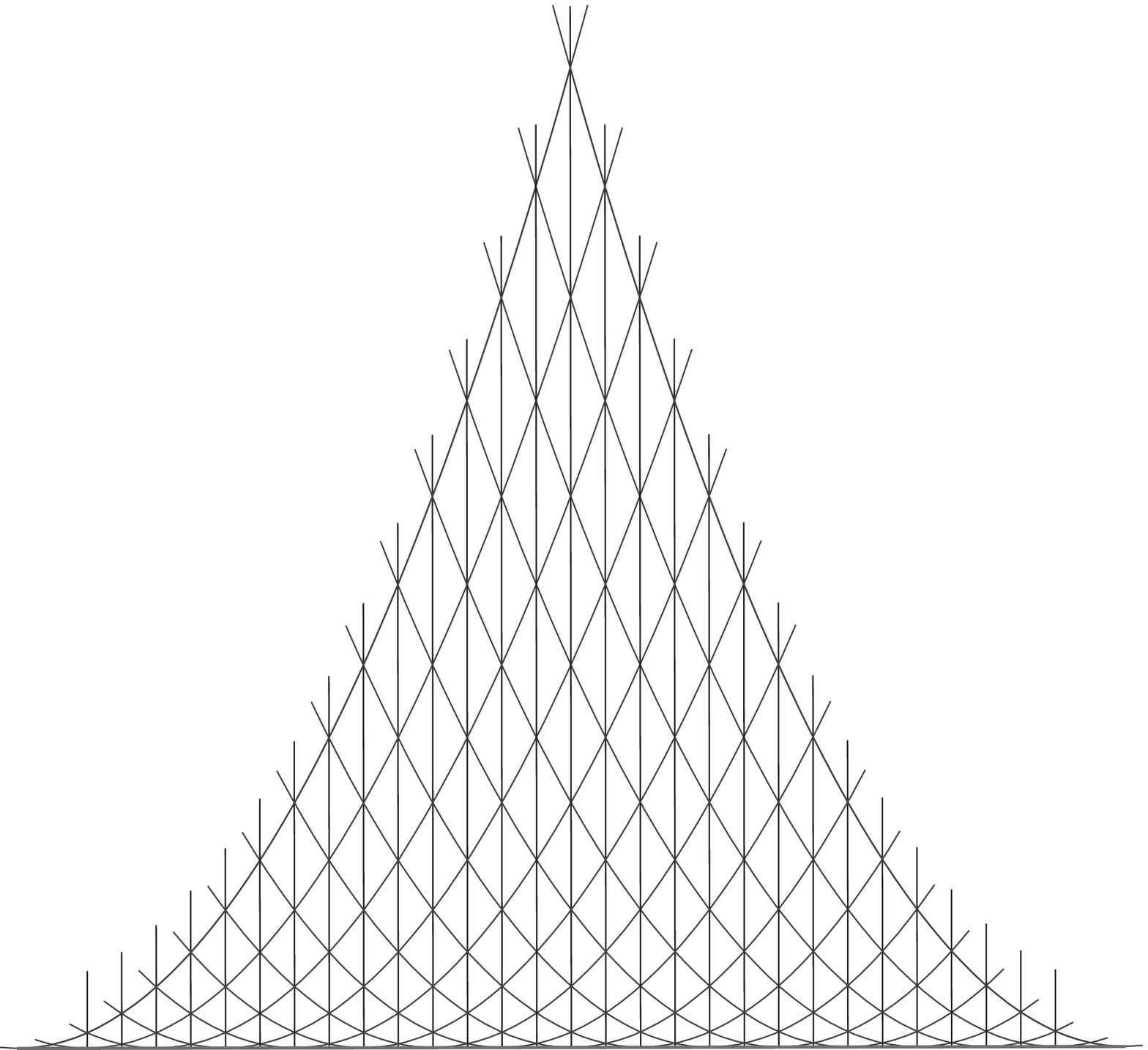,width=48mm} \ \ \ \ \ \ \ \ \ \ \ \ \ \
\ \ \ \ \ \ \ \ \ \epsfig{file=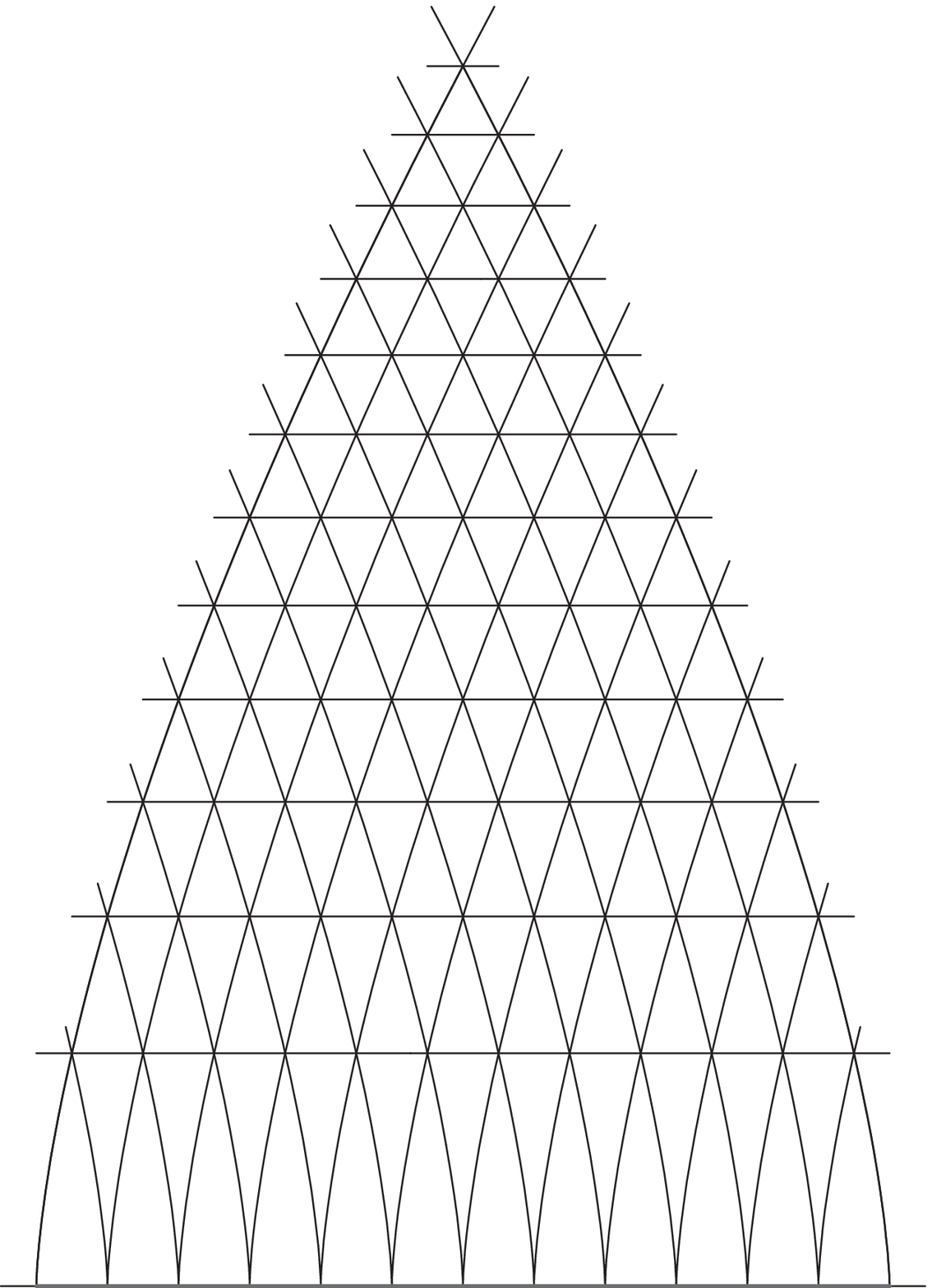,width=32mm}
\caption{Solutions of $p^2=y$ and the lines $dx=0$ (left).
Solutions of $p^2=x$ and the lines $dx=0$ with horizontal y-axis
(right).} \label{Pic_sol_q}
\end{center}
\end{figure}

 \begin{example} {\rm Suppose the criminant of ODE is not Legendrian then this ODE
 is locally equivalent to
\begin{equation}\label{q2}
p^2=x.
\end{equation}
Solutions of this ODE together with the lines $dx=0$ form a
hexagonal web (See Fig. \ref{Pic_sol_q}). The lines $dy=0$  also
complete the 2-web of solutions of (\ref{q2}) to a hexagonal
3-web, but again the surface $M$ of the corresponding cubic
equation $p(p^2-x)=0$ is not smooth at $m=(0,0,0)$.
 }\end{example}

\begin{example} {\rm For completeness let us mention  the case of a regular
point of an implicit cubic ODE. If its 3-web of solutions is
hexagonal it can be mapped to the web of 3 families of parallel
lines $dx=0$, $dy=0$ and $dx+dy=0$. This gives the following
"cubic" ODE
\begin{equation}\label{3l}
p(p+1)=0.
\end{equation}
 }\end{example}
Now we can formulate our classification theorem.
\begin{theorem}\label{classification_th}
Suppose functions $a,b,c$ are real analytic and the following
conditions hold for an implicit cubic ODE
$$
F(x,y,p):=p^3+a(x,y)p^2+b(x,y)p+c(x,y)=0
$$
at a point $m=(x_0,y_0,p_0)\in M :=\{(x,y,p)\in {\mathbb R ^2
\times \mathbb P}^1(\mathbb R): \
F(x,y,p)=0\} $: \\
1) this equation has a hexagonal 3-web of solutions,\\
2) $dF|_{m}\ne 0$,\\
3) ${\rm rank}((x,y,p)\mapsto (F,F_p))|_m=2$ if  $m$ lies on the criminant $C$.\\
Then this ODE is equivalent to one of the following five forms
with respect to some local real analytic isomorphism:
\begin{equation}
\begin{array}{ll}\label{normal}
i) \ \ \ p^3+2xp+y=0, & \mbox{if $p_0$ is a triple root and the criminant is transverse to the} \\
& \mbox{ contact plane field in some punctured neighborhood of $m$,} \\
 ii)\ \ \ p^3+px-y=0, & \mbox{if $p_0$ is a triple root and the criminant is Legendrian,} \\
 iii)\ \ p^2=y, & \mbox{if $p_0$ is a double root and the criminant is Legendrian,} \\
  iv)\ \ p^2=x, & \mbox{if $p_0$ is a double root and the criminant is transverse} \\
& \mbox{to the contact plane at $m$,}\\
  v)\ \ \ p(p+1)=0, & \mbox{if the roots are pairwise distinct at $\pi(m)=(x_0,y_0)$.} \\
\end{array}
\end{equation}
If the functions $a,b,c$ are smooth and conditions 1),2),3) are
satisfied then there is a diffemorphism of a neighborhood of the
point $(x_0,y_0)$ onto a neighborhood of the point $(0,0)$
reducing the above cubic ODE either to one of the four equations
(\ref{normal}ii)-(\ref{normal}v) or to an equations that coincides
with (\ref{normal}i) within the domain, where (\ref{normal}i) has
three real roots.
\end{theorem}
  The main difficulty in proving the above classification theorem brings
the case of irreducible cubic ODE. The idea is to lift its 3-web
of solutions to $M$ and then to the plane $E:\ p_1+p_2+p_3=0$ in
the space of roots of the cubic equation $p^3+A(x,y)p+B(x,y)=0$.
(Note that the general case reduces to this cubic.) Then this
3-web at the plane $E$ has $\mathbb D_3$-symmetry permuting the
roots. Using the regularity condition we construct a $\mathbb
D_3$-equivariant diffeomorphism "upstairs", matching the web to
that of a corresponding normal form. Due to the $\mathbb
D_3$-symmetry the constructed diffeomorphism is lowerable to some
diffeomorphism "downstairs", i.e. to a point transformation in the
plane of solutions. Most of the claims and the proofs below are
given for the smooth case and for some neighborhood of the
projection $\pi(m)$ of $m\in M$, if it is not stated  explicitly.
In section \ref{proof} we discuss how to get rid of the annoying
stipulation in Theorem \ref{classification_th} for the smooth case
(\ref{normal}i) by replacing Definition \ref{def_hex} with a less
geometric one.

\section{Normal forms for a fold point}

In this section we establish normal forms for the case, when the
projection $\pi$ has a fold point at $m$. If cubic equation
(\ref{general_cub}) has two coinciding roots then the the third
root defines a regular point of the projection $\pi$ and the
equation factors to a quadratic equation and a linear one.
Regularity condition (\ref{regularity}) for a double root $p_0$
implies immediately that the projection $\pi$ has a fold point at
$m$. First we find a normal form for fold points and the
symmetries of this normal form. Further we look for the linear in
$p$ (i.e. explicit) equation whose solutions complete the 2-web of
solutions of the quadratic normal form to a hexagonal 3-web.
Finally, we bring these linear terms to some normal forms using
the symmetries of the quadratic equation. We start with a
Legendrian criminant, then consider non-Legendrian criminant and
finally show that the case of an isolated point of tangency of the
criminant and the contact plane is excluded by the regularity
conditions.

\subsection{The case of Legendrian criminant}

\begin{proposition}\label{extension} Consider
 implicit ODE (\ref{implicit}) with a smooth Legendrian criminant and a smooth surface (\ref{surface}).
Then characteristic field (\ref{direction_field}) can be smoothly
extended to the criminant. Moreover, the extended characteristic
field is transverse to the criminant.
\end{proposition}
{\it Proof:} Let $m$ be a point on the criminant. A suitable
contactomorphism $\varphi $ maps $M$ to $M'$ with the following
properties: \\a) the criminant of $M$ is mapped to the line
$x=y=0$,\\b) $\varphi(m)=(0,0,0)$, \\c) the tangent plane $T_mM$
is mapped to the plane $y=0$. \\It suffices to prove the
proposition for the transformed surface $M':=\{(x,y,p)\in \mathbb
R^2 \times \mathbb P^1(\mathbb R): G(x,y,p)=0\}$ with the
Legendrian line $x=y=0$.  Condition c) allows to rewrite
$G(x,y,p)=0$ as
$$
y=g(x,p),$$ while condition a) implies $g(x,p)=xf(x,p)$. Now
$g_p(0,0)=g_x(0,0)=0$ implies $f(0,0)=0$, so that
$f(x,p)=pu(x,p)+xv(x,p)$ by the Hadamard lemma. As the line
$x=y=0$ is Legendrian, the form $dy-pdx$ must vanish on it: $\left
.\{d(xpu(x,p)+x^2v(x,p))-pdx\}\right |_{x=0}=\left
.\{(pu(x,p)dx-pdx\}\right |_{x=0}=0$ or $u(0,p)\equiv 1$. Again by
the Hadamard lemma one gets $u(x,p)=1+xw(x,p)$ and
$$
y=xp+x^2h(x,p)
$$ with $h(x,p)=pw(x,p)+v(x,p)$. Now the characteristic field $\tau$
is defined by restriction of equation $dy-pdx=0$ to $M'$, i.e. by
$\left .\{d(xp+x^2h(x,p))-pdx\}\right
|_{M'}=x\{(1+xh_p(x,p))dp+(2h(x,p)+xh_x(x,p))dx\}=0$. This implies
that the characteristic field on $M'$ in coordinates $(x,p)$ is
generated by the vector field
$(1+xh_p(x,p))\partial_x-(2h(x,p)+xh_x(x,p))\partial_p$ which is
clearly smooth and transverse to the line $x=y=0$ on
$M'$.\vspace{-3ex}\begin{flushright}$\Box$\end{flushright}

\begin{theorem}\label{normal_legendrian}
Let (\ref{implicit}) be an implicit ODE such that corresponding
surface (\ref{surface}) is smooth.  Suppose its criminant $C$ is a
smooth Legendrian curve and the projection $\pi : (x,y,p)\to
(x,y)$ has a fold singularity  at $m\in C$. Then (\ref{implicit})
is locally equivalent to
\begin{equation}\label{normal_legendrian_eq}
p^2=y
\end{equation} with respect to some diffeomorphism $\varphi:
\tilde{U}\to U$, where $\tilde{U},U\in \mathbb R^2$ are
neighborhoods of $(0,0)$, $\pi(m)$ and $\varphi(0,0)=\pi (m)$.
\end{theorem}
{\it Proof:} Due to Proposition \ref{extension} there is a smooth
direction field $\tau$ on $M$. Its integral curves define a
foliation $\mathcal{F}$ of a neighborhood of $m$.  As the
projection $\pi$ has a fold on the criminant, the discriminant
curve $\Delta = \pi(C) \in \mathbb R^2$ is smooth. Let us choose
new coordinates on $(U,\pi(m))$ such that the discriminant curve
turns to the line $y=0$. Then equation (\ref{implicit}) is
equivalent to $y=g(x,p)$. The discriminant curve is a solution
therefore $g(x,p)=pf(x,p)$ by the Hadamard lemma. The criminant
$y=p=0$ is Legendrian hence $ \{d(pf(x,p))-pdx\}_{p=0}=f(x,0)dp=0$
or $f(x,0)\equiv 0$. Now $f(x,p)=ph(x,p)$ and $y=p^2h(x,p)$. As
the projection $\pi$ has a fold at $m\in M$ holds true $ \left.
{\partial _p^2} (p^2h(x,p)) \right |_{x=p=0}=2h(0,0)\ne 0$.
Consequently characteristic field (\ref{direction_field}) on $M$
with $F(x,y,p)=p^2h(x,p)-y$ is generated by the vector field
$$(2ph(x,p)+p^2h_p(x,p))\partial_x-(p^2h_x(x,p)-p)\partial_p=p\{(2h(x,p)+ph_p(x,p))\partial_x+(1-ph_x(x,p))\partial_p\}.$$
Due to $h(0,0)\ne 0$ it is transverse to the kernel of $d\pi$ on
$M$ hence the projection  of each integral curve crossing the
criminant $C$ is smooth and tangent to the discriminant curve.
Locally the equation $ p^2h(x,p)-y=0 $
 can be rewritten as a quadratic
equation
$$
p^2+a(x,y)p+b(x,y)=0.
$$
This easily follows from the Division Theorem since $h(0,0)\ne 0$.
Moreover, as the criminant  $y=p=0$ is a Legendrian curve holds
true $b(x,0)=a(x,0)=0$. Thus one gets $a(x,y)=y\alpha (x,y)$ and
$b(x,y)=y\beta (x,y)$ with $\beta (0,0)\ne 0$ as $M$ is smooth at
$m=(0,0,0)$. Consider $(x,p)$ as local coordinates on $M$ and
define the map $i:M\to M$ by
\begin{equation}\label{involution}
i(x,p)=(x,-p-p^2h(x,p)\alpha (x, p^2h(x,p))).
\end{equation}
This map is the involution that permutes the roots $p_1,p_2$ of
our quadratic ODE. Let the foliation $\mathcal{F}$ on $M$ be
defined locally by $I(x,p)=const$ with $\left. {\rm grad}
(I)\right |_{x=p=0}\ne 0$, where $I$ is a first integral of the
characteristic field $\tau$. Then the functions  $I$ and
$J:=i^{*}(I)$ are functionally independent as $di(\tau) $ and
$\tau$ are transverse to each other. Since $\tau$ is transverse
also to the criminant the partial derivative $\left. I_x\right
|_{x=p=0}$ does not vanish. Hence one can chose $I$ so that
$I(x,0)=x$ which implies $J(x,0)=x$. Let us take thus normalized
functions $I,J$ as local coordinates on $M$. Note that the
following
 relation holds true
\begin{equation}\label{lowerability}
 \pi(I,J)=\pi(J,I),
\end{equation}
  since $i(I,J)=(J,I)$. For the
 ODE $\tilde{p}^2=\tilde{y}$ the above defined
 objects are as follows: $\tilde{I}(\tilde{x},\tilde{p})=\tilde{x}-2\tilde{p}$, $\tilde{i}(\tilde{x},\tilde{p})=(\tilde{x},-\tilde{p})$,
 $\tilde{J}(\tilde{x},\tilde{p})=\tilde{x}+2\tilde{p}$. Now define a diffeomorphism germ $\psi:\tilde{M},0\to
 M,0$ by $\psi
 (\tilde{I},\tilde{J})=(\tilde{I},\tilde{J})$.  By
 (\ref{lowerability}) there exists such a map germ $\varphi:\mathbb R^2,0\to \mathbb
 R^2,0$ that the following diagram commutes:
\begin{equation}
\begin{array}{ccc}
  M,0 & \stackrel{\psi}{\longleftarrow} &\tilde{ M},0\\
  \downarrow \pi &   &\downarrow \tilde{\pi}\\
  \mathbb R^2,0 & \stackrel{\varphi}{\longleftarrow}& \mathbb R^2,0\\
\end{array}%
\end{equation}
We claim that the map germ $\varphi$ is the searched for
diffeomorphism. By construction it maps solutions to solutions. To
show that $\varphi$ is differentiable consider the map germ $\pi
\circ \psi: \tilde{M},0 \to \mathbb R^2,0$. Applying Malgrange's
Preparation Theorem to this map germ in coordinates
$(\tilde{x},\tilde{p})$ on $\tilde{M},0$ and $(x,y)$ on $\mathbb
R^2,0$ one gets
$$
x=X_0(\tilde{x},\tilde{p}^2)+\tilde{p}X_1(\tilde{x},\tilde{p}^2),
\ \
y=Y_0(\tilde{x},\tilde{p}^2)+\tilde{p}Y_1(\tilde{x},\tilde{p}^2),
$$
where the functions $X_0,X_1,Y_0,Y_1$ are smooth. Since $\pi \circ
\psi (\tilde{x},\tilde{p})=\pi \circ
\psi(\tilde{x},-\tilde{p})=(x,y)$ holds true
$X_1(\tilde{x},\tilde{p}^2)=Y_1(\tilde{x},\tilde{p}^2)=0$. With
$\tilde{p}^2=\tilde{y}$ the above formulas take the form
$$
x=X_0(\tilde{x},\tilde{y}), \ \ y=Y_0(\tilde{x},\tilde{y}),
$$
 and therefore define a smooth map germ $\varphi$. Applying the same
 considerations to the map germ $\tilde{\pi}\circ \psi ^{-1}$ we see that $\psi$ has an inverse for $y \ge
0.$ Therefore $\psi$  is invertible.\ \ \ \ \ \ \ \ \ \ \ \ \ \ \ \ \ \ \ \ \ \ \ \ \ \ \ $\Box$\\

\noindent {\bf Remark.} Note that equation
(\ref{normal_legendrian_eq}) has a nontrivial symmetry group. For
example, the scaling $x\to \lambda x, \ \ y\to \lambda ^2 y$
leaves it invariant. Therefore the diffeomorphism $\varphi$ in
Theorem \ref{normal_legendrian} is not unique.

\begin{proposition}\label{symmetry_quadratic_legendrian}
Equation (\ref{normal_legendrian_eq}) has an infinite symmetry
pseudogroup. Its transformations are given by
\begin{equation}\label{finite_symmetry_quadratic_legendrian}
\tilde{x}=F(x+2\sqrt y)+F(x-2\sqrt y), \ \ \tilde{y}=\frac{1}{4} [
F(x+2\sqrt y)-F(x-2\sqrt y) ]^2.
\end{equation}
Here $F$ is a smooth function subjected to $\left. F'(u)\right
|_{u=0}\ne 0$.
 Infinitesimal
generators of this pseudogroup have the form
\begin{equation}\label{infinitesimal_symmetry_quadratic_legendrian}
\{f(x+2\sqrt y)+f(x-2\sqrt y)\}\partial _x+\sqrt y\{f(x+2\sqrt
y)-f(x-2\sqrt y)\}\partial _y,
\end{equation}
where $f$ is an arbitrary smooth function.
\end{proposition}
{\it Proof:} Consider the action of symmetry group transformation
on $M$ in coordinates $(u,v)$, where
\begin{equation}\label{uv}
u=x-2p,\ \ v=x+2p.
\end{equation}
To preserve the foliations $\mathcal{F}$ and $i(\mathcal{F})$
determined by the  direction fields $\tau$ and $di(\tau)$ it must
have the following form $\bar{u}=2F(u),\ \bar{v}=2G(v)$. This
transformation is a symmetry if it commutes with $i$. Note that
$i$ permutes $u$ and $v$ and sends $\partial _p$ to $-\partial_p$.
Hence $F=G$. Now substitutions $\bar{x}=\frac{1}{2}(F(v)+F(u))$
and $\bar{y}=\bar{p}^2=\frac{1}{4}(F(v)+F(u))^2$ gives
(\ref{finite_symmetry_quadratic_legendrian}). Condition $\left.
F'(u)\right |_{u=0}\ne 0$ is equivalent to non-vanishing of the
Jacobian: $\left. \frac{\partial (\tilde{x},\tilde{y})}{\partial
(x,y)}\right |_{x=y=0}\ne 0.$

 Consider now  infinitesimal symmetries of
(\ref{normal_legendrian_eq}). They are defined by operators $\xi
(x,y)\partial _x+\eta (x,y)\partial _y$. These operators must be
liftable to $M$. On $M$ the lifted vector field must be a symmetry
of foliations $\mathcal{F}$ and $i(\mathcal{F})$. In coordinates
$(u,v)$  any infinitesimal symmetry $X$ of foliations
$\mathcal{F}$ and $i(\mathcal{F})$ is easy to write down:
\begin{equation}\label{X_uv}
X=f(u)\partial _u+g(v)\partial _v.
\end{equation}
In coordinates $(x,p)$ on $M$ this operator $X$ takes the form
$$
X=\frac{1}{2}(g(v)+f(u))\partial _x+\frac{1}{4}(g(v)-f(u))\partial
_p.
$$
$X$ is lowerable iif $di(X)=X$.  Now lowerability condition
amounts to  $f(u)+g(v)=f(v)+g(u)$ and $g(v)-f(u)=f(v)-g(u)$. Thus
the first equation gives $f(u)-g(u)=c=const$ and the second
implies $c=0$. Substitution $u=x-2\sqrt y$, $v=x+2\sqrt y$,
$p=\sqrt y$ into $X$ gives
(\ref{infinitesimal_symmetry_quadratic_legendrian}) for $y >0 $.
The obtained transformations are correctly defined for $y> 0$ and
are smoothly (analytically) extendable for $y\ge 0$. The
possibility to extend them for $y\le 0$ easily follows from
Malgranges's Preparation Theorem: one considers the projection
$\pi:M\to \mathbb R^2$, $(x,p)\mapsto (x,p^2)$  in coordinates
$x,p$ on $M$ and $x,y$ on $\mathbb R^2$ and observes that
$F(x+2p)+F(x-2p)$ is an even and $F(x+2p)-F(x-2p)$ is an odd
function with respect to $p$.
\vspace{-4ex}\begin{flushright}$\Box$\end{flushright} \noindent To
use the above symmetries we need the following lemma.
\begin{lemma}\label{normal_1_dim_vector}
If a smooth (analytic) function germ $f: \mathbb R,0\to \mathbb
R,0$ is not flat then the vector field germ $f(u)\partial _u$ on
$\mathbb R,0$ is equivalent to $u^k\partial _u$ with respect to a
certain smooth (analytic) coordinate transformation
$\bar{u}=F(u)$, where $k\in \mathbb N$ is such that
$\left.\frac{d^k f(u)}{{d u}^k}\right |_{u=0}$ is the first
non-vanishing derivative at $0$  .
\end{lemma}
{\it Proof:} $F(u)$ must satisfy ODE $f(u)F'(u)=(F(u))^k$. The
existence of $F(u)$ with $F'(u)|_{u=0}\ne 0$ is easily verified in
both smooth and analytic cases.
\vspace{-4ex}\begin{flushright}$\Box$\end{flushright}
\begin{theorem}\label{quadratic_legendrian_web}
Suppose the solutions of (\ref{normal_legendrian_eq}) and those of
\begin{equation}\label{explicit}
\alpha (x,y)dx+\beta (x,y)dy=0,
\end{equation}
where $\alpha (x,y),\beta (x,y)$ are non-flat functions at
$(0,0)$, form together a hexagonal 3-web. Then there is a local
symmetry of (\ref{normal_legendrian_eq}) at $(0,0)$ that maps
equation (\ref{explicit}) to one of the two following forms for
$y\ge 0$:
\begin{equation}\label{normal_explicit_legendrian}
\begin{array}{c}
  [(x+2\sqrt y)^k+(x-2\sqrt y)^k]dx-\frac{1}{\sqrt y}[(x+2\sqrt
y)^k-(x-2\sqrt y)^k]dy=0,\ {\rm or} \\
 \\
\sqrt y[(x+2\sqrt y)^k-(x-2\sqrt y)^k]dx-[(x+2\sqrt y)^k+(x-2\sqrt
y)^k]dy=0,\\
\end{array}
\end{equation}
where $k$ is a non-negative integer.  In particular, if $\alpha
(x,y),\beta (x,y)$ are non-flat functions with $(\alpha
(0,0),\beta (0,0))\ne (0,0)$,
 one gets three
normal  forms:
\begin{equation}\label{to_quadratic_legendrian}
  a) \ dx=0, \ \  \ b) \ dy=0, \ \ \ c)\ 2dy-xdx=0.
\end{equation}
Moreover, if equation (\ref{explicit}) for $y\ge 0$ is equivalent
to (\ref{to_quadratic_legendrian}a) or
(\ref{to_quadratic_legendrian}b) then it can be reduced to
(\ref{to_quadratic_legendrian}a), respectively
(\ref{to_quadratic_legendrian}b) by a symmetry of
(\ref{normal_legendrian_eq}) in some neighborhood of the point
$(0,0)$
\end{theorem}
{\it Proof:} Let us introduce operators $U,V$ of differentiation
along the curves of the foliations $\mathcal{F}$ and
$i(\mathcal{F})$ on $M$:
\begin{equation}\label{UV}
U=\partial _p+2\partial _x, \ \ V=\partial _p-2\partial _x.
\end{equation}
Then these operators commute and satisfy the following relations:
$$
U(u)=0,\ \  U(v)=4, \ \  V(u)=-4,\ \ V(v)=0.
$$
Consequently a  direction field on $M$, whose integral curves form
a hexagonal 3-web together with $\mathcal{F}$ and
$i(\mathcal{F})$, must be generated by a vector field, commuting
with $U$ and $V$ (for the details see \cite{BB}, p.17). Such a
vector field has the form
$$
Y=f(v)U+g(u)V.
$$
The  direction field generated by $Y$ is the lift to $M$ of the
 direction field induced by (\ref{explicit}) iff $Y\wedge di(Y)=0$. This
gives $g=\pm f$. Projecting from $M$ to the plane  one obtains
$$
\frac{1}{\sqrt y}[f(x+2\sqrt y)-f(x-2\sqrt
y)]\partial_x+[f(x+2\sqrt y)+f(x-2\sqrt y)]\partial_y
$$
for $g=f$ and
$$
[f(x+2\sqrt y)+f(x-2\sqrt y)]\partial _x+\sqrt y[f(x+2\sqrt
y)-f(x-2\sqrt y)]\partial_y,
$$
for $g=-f$. Now applying symmetry
(\ref{finite_symmetry_quadratic_legendrian}) with $F$ satisfying
$f(u)F'(u)=(F(u))^k$ (see  Lemma \ref{normal_1_dim_vector}) we
reduce (\ref{explicit}) to one of the forms
(\ref{normal_explicit_legendrian}) for $y\ge 0$. If the found
symmetry maps (\ref{explicit}) to (\ref{to_quadratic_legendrian}a)
than we can construct a diffeomorphism $\varphi$ such that it is
the identity for $y\ge 0$ and maps integral curves of
(\ref{explicit}) to the lines $x=const$ as follows. As equation
(\ref{explicit}) is not singular at $(0,0)$ it has a smooth first
integral $I(x,y)$ coinciding with $x$ for $y\ge 0$. Then $\phi$ is
defined by $(x,y)\mapsto (I(x,y),y)$.   Similarly, for
(\ref{to_quadratic_legendrian}b) we define $\phi$  by
$(x,y)\mapsto (x,I(x,y)$, where $I(x,y)$ is the first integral of
(\ref{explicit}), coinciding with $y$ for $y\ge 0$.
\vspace{-3ex}\begin{flushright}$\Box$\end{flushright} \noindent

\subsection{The case of non-Legendrian criminant}

\begin{theorem}\label{normal_non_legendrian}
Let (\ref{implicit}) be an implicit ODE such that the
corresponding surface $M$  is smooth.  Suppose its criminant $C$
is a smooth curve, the projection $\pi : (x,y,p)\to (x,y)$ has a
fold singularity  at $m\in C$, and the contact plane is not
tangent to $M$ at $m$. Then (\ref{implicit}) is locally equivalent
to
\begin{equation}\label{normal_non_legendrian_eq}
p^2=x
\end{equation} with respect to some diffeomorphism $\varphi:
\tilde{U}\to U$, where $\tilde{U},U\in \mathbb R^2$ are
neighborhoods of $(0,0)$, $\pi(m)$ and $\varphi(0,0)=\pi (m)$.
\end{theorem}
The proof is given in \cite{Ag}, p.27. Similar to  the case of
Legendrian criminant, the diffeomorphism $\varphi$  is not unique.
\begin{proposition}\label{symmetry_quadratic_non_legendrian}
Equation (\ref{normal_non_legendrian_eq}) has an infinite symmetry
pseudogroup. Its transformations are given by
\begin{equation}\label{finite_symmetry_quadratic_non_legendrian}
\tilde{x}=\sqrt[3]{\frac{1}{16}[ F(3y+2x\sqrt x)-F(3y-2x\sqrt x)]
^2}, \ \ \tilde{y}=\frac{1}{6} (F(3y+2x\sqrt x)+F(3y-2x\sqrt x)).
\end{equation}
Here $F$ is a smooth function subjected to $\left. F'(u)\right
|_{u=0}\ne 0$.
 Infinitesimal
generators of this pseudogroup have the form
\begin{equation}\label{infinitesimal_symmetry_quadratic_non_legendrian}
\frac{1}{\sqrt x}\{ f(3y+2x\sqrt x)-f(3y-2x\sqrt x)\}
 \partial _x+\{f(3y+2x\sqrt x)+f(3y-2x\sqrt x)\}\partial _y,
\end{equation}
where $f$ is an arbitrary smooth function.
\end{proposition}
{\it Proof:} On the surface $M$ defined by (\ref{surface}) we
choose $(y,p)$ as local coordinates. Solutions of
(\ref{normal_non_legendrian_eq}) define foliations $\mathcal{F}$
and $i(\mathcal{F})$  by
\begin{equation}\label{uv1}
u:=3y-2p^3=const, \ \ \ v:=3y+2p^3=const,
\end{equation}
where $i$ is the involution $i: M\to M,\ (x,y,p)\to (x,y,-p)$.

To prove the finite transformation formulas observe that the
symmetry group transformation $(x,y)\mapsto
(\tilde{x},\tilde{y})$ lifted to $M$ must satisfy
$$
3\tilde{y}+2\tilde{p}^3=F(3y+2p^3), \ \ \
3\tilde{y}-2\tilde{p}^3=F(3y-2p^3).
$$
In fact, to preserve the foliations $\mathcal{F}$ and
$i(\mathcal{F})$ it is necessary that
$3\tilde{y}+2\tilde{p}^3=F(3y+2p^3), \ \
3\tilde{y}-2\tilde{p}^3=G(3y-2p^3).$ On the line $p=0$ one has
$F(3y)=G(3y)=3\tilde{y}$. This implies
(\ref{finite_symmetry_quadratic_non_legendrian}). The condition
$\left. F'(u)\right |_{u=0}\ne 0$ is equivalent to non-degeneracy
of the Jacobian of the transformation.

Consider now an infinitesimal symmetry $\xi (x,y)\partial _x+\eta
(x,y)\partial _y$. Lifted on $M$ it turns to
$$X=\eta (p^2,y)\partial _y+g(y,p)\partial _p.$$ (On can write down
an explicit expression for $g$, but we do not need it.) Since $X$
is a symmetry of (\ref{normal_non_legendrian_eq}) it must satisfy
$$X(u)=3f(u)=3f(3y-2p^3), \ \ \ X(v)=3k(v)=3k(3y+2p^3),$$ for some smooth functions $f,k$.  This is equivalent to
$$
\eta (p^2,y)=2p^2g(y,p)+f(3y-2p^3), \ \ \  \eta
(p^2,y)=-2p^2g(y,p)+k(3y+2p^3).
$$
Substituting $p=0$ into the difference of the above equations
$$
4p^2g(y,p)=k(3y+2p^3)-f(3y-2p^3)$$ one gets $k=f$. Hence the
functions $\eta $ and $g$  are well defined by
$$
\eta(p^2,y)=\frac{1}{2}(f(3y+2p^3)+f(3y-2p^3)),\ \ \
g(y,p)=\frac{1}{4p^2}(f(3y+2p^3)-f(3y-2p^3)).
$$
Note that $g(y,0)=0$, hence  $X$ is tangent to the criminant and
therefore lowerable (see \cite{Aw}). Up to scaling the lowered
operator $X$  becomes
(\ref{infinitesimal_symmetry_quadratic_non_legendrian}). The
extension of the defined transformation for $x\le 0$ is again
justified by Malgrange's Preparation Theorem. (See the detail in
the proof of Theorem \ref{symmetry_quadratic_legendrian}.)
\vspace{-3ex}\begin{flushright}$\Box$\end{flushright}
\begin{theorem}\label{quadratic non_legendrian_web}
Suppose the solutions of (\ref{normal_non_legendrian_eq}) and
those of
\begin{equation}\label{explicit_1}
\alpha (x,y)dx+\beta (x,y)dy,
\end{equation}
where $\alpha (x,y),\beta (x,y)$ are non-flat functions  at
$(0,0)$, form together a hexagonal 3-web. Then there is a local
symmetry of (\ref{normal_non_legendrian_eq}) that maps equation
(\ref{explicit_1}) to one of the two following forms for $x\ge 0$:
$$
\frac{1}{\sqrt x}{[(3y+2x\sqrt x)^k-(3y-2x\sqrt x)^k]}dx
-\sqrt[3]{16}[(3y+2x\sqrt x)^k+(3y-2x\sqrt x)^k]dy=0\ {\rm or}
$$
$$
[(3y+2x\sqrt x)^k+(3y-2x\sqrt x)^k]dx-\sqrt[3]{16}\sqrt x
[(3y+2x\sqrt x)^k-(3y-2x\sqrt x)^k]dy=0,
$$
where $k$ is non-negative integer.  In particular, if $\alpha
(x,y),\beta (x,y)$ are non-flat functions with $(\alpha
(0,0),\beta (0,0))\ne (0,0)$, then equation (\ref{explicit_1}) can
be reduced to one of the following two normal forms in some
neighborhood of the point $(0,0)$:
\begin{equation}\label{to_quadratic_non_legendrian}
  a) \ dx=0, \ \  \ b) \ dy=0.
\end{equation}
\end{theorem}
{\it Proof:} Let us introduce operators $U,V$ of differentiation
along the curves of the foliations $\mathcal{F}$ and
$i(\mathcal{F})$:
\begin{equation}\label{UV}
U=\partial _p+2p^2\partial _y, \ \ V=\partial _p-2p^2\partial _y.
\end{equation}
Then the operators $\frac{1}{p^2}U$ and $\frac{1}{p^2}V$ commute
and satisfy the following relations:
$$
\frac{1}{p^2}U(u)=0,\ \ \frac{1}{p^2}U(v)=12, \ \
\frac{1}{p^2}V(u)=-12,\ \ \frac{1}{p^2}V(v)=0.
$$
A  direction field on $M$, whose integral curves form a hexagonal
3-web together with $\mathcal{F}$ and $i(\mathcal{F})$, must be
generated by the vector field, commuting with $\frac{1}{p^2}U$ and
$\frac{1}{p^2}V$. Such a vector field has the form
$$
Y=f(v)\frac{1}{p^2}U+g(u)\frac{1}{p^2}V.
$$
The  direction field generated by $Y$ is the lift to $M$ of the
 direction field induced by (\ref{explicit_1}) iff $Y\wedge di(Y)=0$.
This gives $g=\pm f$ (compare with the proof of Theorem
\ref{quadratic_legendrian_web}). Projecting from $M$ to
$(x,y)$-plane  one obtains
$$
[f(3y+2x\sqrt x)+f(3y-2x\sqrt x)]\partial_x+\sqrt x[v(3y+2x\sqrt
x)-f(3y-2x\sqrt x)]\partial_y
$$
for $g=f$ and
$$
\frac{1}{\sqrt x}[f(3y+2x\sqrt x)-f(3y-2x\sqrt
x)]\partial_x+[v(3y+2x\sqrt x)+f(3y-2x\sqrt x)]\partial_y
$$
for for $g=-f$. Now applying symmetry
(\ref{finite_symmetry_quadratic_non_legendrian}) with $F$
satisfying $f(u)F'(u)=(F(u))^k$ (see  Lemma
\ref{normal_1_dim_vector}) we complete the proof. The details can
be found in the proof of Theorem \ref{quadratic_legendrian_web}.
\vspace{-3ex}\begin{flushright}$\Box$\end{flushright} {\bf Remark
1.} Real analytic versions of Theorems
\ref{quadratic_legendrian_web} and \ref{quadratic
non_legendrian_web} are true without the stipulations $y\ge 0$ and
$x\ge 0$ respectively.\\
{\bf Remark 2.} One can not extend the claim of Theorem
\ref{quadratic_legendrian_web}  for $y<0$ in the smooth  case for
equation (\ref{to_quadratic_legendrian}c). Its solutions are
parabolas $y=x^2/4+C$. They cross the line $y=0$ in two points if
$C<0$. If one smoothly deforms equation
(\ref{to_quadratic_legendrian}) in the domain $y<0$ then the
solution of the deformed equation starting from some point
$(x_0,0)$ with $x_0<0$ will not necessarily pass through
$(-x_0,0)$, i.e. this solution returns to a "wrong parabola".
\begin{corollary}\label{no_isolated_tangency}
If the following conditions hold for implicit ODE (\ref{general_cub}) at a point $m=(x_0,y_0,p_0)\in M$: \\
1) ODE (\ref{implicit}) has a hexagonal 3-web of solutions,\\
2) $p_0$ is a double root of (\ref{implicit}) at $\pi(m)=(x_0,y_0)$, \\
3) regularity condition (\ref{regularity}) is satisfied, i.e. ${\rm rank}((x,y,p)\mapsto (F,F_p))|_m=2$,\\
then its criminant is either Legendrian or transverse to the
contact plane field in some neighborhood of $m$
\end{corollary}
{\it Proof:} Denote by $C_t$ the closed set of points on the
criminant $C$, where the contact plane is tangent to $C$. Suppose
 $m$ is not a point of $C\setminus C_t$ and not an interior point of
 $C_t$. Then $m$ is a boundary point of $C_t$.
Now Theorem \ref{quadratic non_legendrian_web} implies that for
each point $m'$ sufficiently close to $m$ and such that $m'\ne m$,
$m'\notin C_t$ equation (\ref{implicit}) is locally equivalent to
a product of an explicit ODE $dx=0$  and quadratic equation
(\ref{normal_non_legendrian_eq}), i.e. the solutions of the linear
factor are tangent to the discriminant curve at $\pi(m')$ and
therefore at $\pi(m)$. Further Theorem
\ref{quadratic_legendrian_web} implies that for each point $m'$
sufficiently close to $m$ and  such that $m'\ne m$, $m'\in C_t$
equation (\ref{implicit}) is locally equivalent to a product of an
explicit ODE $dx=0$  and quadratic equation
(\ref{normal_legendrian_eq}), i.e. the solutions of the quadratic
factor are tangent to the discriminant curve at $\pi(m')$ and
therefore at $\pi(m)$. But that means that the root $p_0$ is
triple.  Thus our assumption is false and the corollary is proved.
\ \ $\Box$\\

\noindent{\bf Remark 3.} As follows from the above proof the
hypothesis of Corollary \ref{no_isolated_tangency} also implies
that there is no isolated points of tangency of the contact plane
and the criminant.

\section{Normal form for an ordinary cusp point}
In this section we use the results of the previous one to
establish normal forms for the case of a cusp singularity of the
projection $\pi$ on $M.$ Regularity condition (\ref{regularity})
for a triple root $p_0$ implies immediately that the projection
$\pi$ has a cusp point at $m=(x_0,y_0,p_0)\in M$. We start with a
Legendrian criminant, then consider non-Legendrian criminant and
finally show that one can not "glue"  Legendrian criminant with
non-Legendrian one at the cusp point.

\begin{lemma}\label{cubic}
If the following conditions hold for implicit ODE (\ref{implicit}) at a point $m=(x_0,y_0,p_0)\in M$: \\
1) ODE (\ref{implicit}) has a hexagonal 3-web of solutions,\\
2) $p_0$ is the triple root of (\ref{implicit}) at $\pi(m)=(x_0,y_0)$, \\
3) regularity condition (\ref{regularity}) is satisfied, i.e. ${\rm rank}((x,y,p)\mapsto (F,F_p))|_m=2$,\\
then it is locally equivalent to
\begin{equation}\label{CUB}
p^3+A(x,y)p+B(x,y)=0,
\end{equation}
where \\
1)  the projection $\pi$ has an ordinary cusp singularity at
$(0,0,0)$ with $A(0,0)=B(0,0)=0$,\\
2) $A,B$ are local coordinates at $(0,0)$, i.e.
$\left. \frac{\partial(A,B)}{\partial(x,y)}\right |_{x=y=0}\ne 0$\\
3) $B_x(0,0)=0$.
\end{lemma}
{\it Proof:} Since equation (\ref{implicit}) has the triple root
$p_0$ at $\pi(m)$ it is locally equivalent to some cubic equation
(\ref{general_cub}). Further, the coefficient by $p^2$ in this
cubic equation is  killed by a coordinate transform of the form
$y=f(\tilde{x},\tilde{y})$, $x=\tilde{x}$, satisfying
$$
3f_x(x,y)+a(x,y)=0.
$$
This transform respects the regularity conditions.  Thus our
implicit equation $F=0$ becomes
$$
F(x,y,p)=p^3+A(x,y)p+B(x,y)=0.
$$
 Without loss of generality it
can be assumed that $(x_0,y_0)=(0,0)$. As the equation
$p_0^3+A(0,0)p_0+B(0,0)=0$ has a triple root holds $p_0=0$.
Therefore the functions $A,B$ must also vanish at $(0,0)$. Now
regularity condition (\ref{regularity}) at $m=(0,0,0)\in C$ reads
as
$$
{\rm rank}\left.\left(%
\begin{array}{ccc}
  A_xp+B_x & A_yp+By & 3p^2+A \\
  A_x & A_y & 6p \\
\end{array}%
\right)\right |_{x=y=p=0} =\left.\left(%
\begin{array}{ccc}
  B_x & By & 0\\
  A_x & A_y & 0 \\
\end{array}
\right)\right |_{x=y=p=0}=2.
$$
Thus claims 1) and 2) are proved. Moreover, the discriminant curve
$\Delta=\pi (\{(x,y,p):\ p^3+A(x,y)p+B(x,y)=3p^2+A(x,y)=0 \})$ has
an ordinary cusp at $\pi(m)=(x_0,y_0)$.

If solutions of equation (\ref{CUB}) form a hexagonal 3-web the
curvature of this 3-web must vanish identically. This is
equivalent to the following cumbersome partial differential
equation for the functions $A,B$:
\begin{equation}\label{PDE}
\begin{array}{l}
  (4A^3+27B^2)(9BA_{xx}-2A^2A_{xy}+6ABA_{yy}-6AB_{xx}-9BB_{xy}-4A^2B_{yy})+\\
  +108A^2BA_xB_y-108AB^2A_xA_y+162B^3A_y^2+40A^4A_yB_y-108A^2BA_x^2+\\
  +216A^2BB_y^2-36A^3B_xB_y+108A^2BA_yB_x-378AB^2A_yB_y-405B^2A_xB_x+\\
  -48A^3BA_y^2+8A^4A_xA_y+243B^2B_xB_y+84A^3A_xB_x+324ABB_x^2=0.\\
\end{array}
\end{equation}
This equation is obtained by direct lengthy but straightforward
computation. (Expressions for the corresponding web curvature for
a cubic ODE can be also found in \cite{Mr} and \cite{Nw}). As was
shown above the functions $A,B$ can be taken as local coordinates
around $(0,0)$. Then all partial derivatives of $A,B$ with respect
to $x$ and $y$ are smooth functions of $A,B$. The homogeneous part
of second order of Taylor expansion of l.h.s. of (\ref{PDE})
around $(0,0)$ is
$$-405B^2A_x(0,0)B_x(0,0)+243B^2B_x(0,0)B_y(0,0)+324AB(B_x(0,0))^2.$$
It must vanish. In particular, $B_x(0,0)^2=0$ as the coefficient
by $AB$. \vspace{-3ex}\begin{flushright}$\Box$\end{flushright}

\subsection{The case of Legendrian criminant}

\begin{theorem}\label{Clairaut_th}
If the following conditions hold for implicit ODE (\ref{implicit}) at a point $m=(x_0,y_0,p_0)\in M$: \\
1) ODE (\ref{implicit}) has a hexagonal 3-web of solutions,\\
2) $p_0$ is the triple root of (\ref{implicit}) at $\pi(m)=(x_0,y_0)$,\\
3) its criminant $C$ is a Legendrian curve with ${\rm rank}((x,y,p)\mapsto (F,F_p))|_C=2$, \\
then it is locally equivalent to the following Clairaut equation
\begin{equation}\label{Clairaut}
P^3+PX-Y=0.
\end{equation}
\end{theorem}
{\it Proof:} Lemma \ref{cubic} reduces equation (\ref{implicit})
to  (\ref{CUB}) with $A(0,0)=B(0,0)=B_x(0,0)=0$ and $\left.
\frac{\partial(A,B)}{\partial(x,y)}\right |_{x=y=0}\ne 0$. Thus
the tangent plane $T_mM$ to the surface $M$ at $m=(0,0,0)$ is the
plane $y=0$ and the tangent line to the discriminant curve at
$(0,0)$ is $y=0$.  (The condition $B_x(0,0)=0$ for the case of
Legendrian criminant can be obtain also by the following
geometrical consideration. By Proposition \ref{extension} the
characteristic field $\tau$, given by (\ref{direction_field}), can
be smoothly extended to the criminant $C$. As $\tau$ is transverse
to $C$ the projection of integral curves of $\tau$  are smooth
curves in $\mathbb R^2$ tangent to the discriminant curve $\Delta$
at the origin $(0,0)$. As $p_0=0$ this implies the claim.) With
$B_x(0,0)=0$ one gets from the regularity condition
$$
A_x(0,0)\ne 0,\ \ B_y(0,0)\ne 0.
$$
Therefore one can choose $p,A$ as local coordinates on the surface
$M$ at $m=(0,0,0)$ and $A,B$ as local coordinates on $\mathbb
R^2,0$. In these coordinates the projection $\pi$ is the Whitney
map. The criminant is parameterized by $p$ as follows
$$
A=-3p^2, \ \ B=2p^3.
$$
Its projection is the discriminant curve
$\Delta:=\{(A,B):27B^2+4A^3=0\}$. The set of points projected to
the discriminant curve is the criminant itself and the following
curve
\begin{equation}\label{preimage_discriminant}
D:= \{(p,A)\in M:A=-\frac{3}{4}p^2\}.
\end{equation}
This follows from the observation that the value $-2p$ is the
third root of (\ref{CUB}) at the discriminant, where the double
root is $p.$ The curve $D$ is tangent to $C$ at $0$, thus the
characteristic field $\tau$ is transverse also to $D$.

Consider the following map
\begin{equation}\label{pull_to_pq}
f:(\mathbb R^2,0)\to (M,0),\ (p,q)\mapsto  (p,A)= (p,
-q^2-pq-p^2).
\end{equation}
This map has a fold singularity on the line $L_1:=\{(p,q):\ p+2q=0
\}$. This line is mapped by $f$ to the curve $D$ since
$-q^2-pq-p^2=-\left (q+\frac{p}{2}\right )^2-\frac{3}{4}p^2$. Note
that if $f(p,q)=(p,A)$ then  $f^{-1}(p,A)=\{(p,q)\cup (p,-p-q)\}.$
The pull back $\tilde{\tau}$ of the characteristic field $\tau$ by
$f^*$ from $M,0$ to $\mathbb R^2,0$ must be tangent to the kernel
of $df$, i.e. to the vector field $\partial _q$, since $\tau$ is
transverse to $D$. Moreover, the foliation $\mathcal{F}_1$ of
integral curves of $\tilde{\tau}$ is invariant with respect to the
following linear involution
$$
g_1:\mathbb R^2 \to \mathbb R^2,\ (p,q)\mapsto (p,-p-q).
$$
Consider also the following two linear involutions
\begin{equation}\label{g2_g3}
g_2:\mathbb R^2 \to \mathbb R^2,\ (p,q)\mapsto (-p-q,q),\ \ \
g_3:\mathbb R^2 \to \mathbb R^2,\ (p,q)\mapsto (q,p).
\end{equation}
The linear maps $g_1,g_2,g_3$ generate the group $\mathbb D_3$,
the symmetry group of equilateral triangle, which can be viewed as
the group of linear transformations of the plane $p+q+r=0$ in
$\mathbb R ^3$, generated by permutations of the coordinates
$(p,q,r)$ in $\mathbb R ^3$. The orbit of a point $(p,q)$ under
this group action is the inverse image of the point
$(A,B)=\pi(f(p,q))\in \mathbb R^2$ under the Vieta map $V:=\pi
\circ f$. Therefore the three foliations $\mathcal{F}_1$,
$\mathcal{F}_2:=g_2(\mathcal{F}_1)$ and
$\mathcal{F}_3:=g_3(\mathcal{F}_1)$ form a hexagonal 3-web.
Moreover, this 3-web is not singular at $(0,0)$ and has the
symmetry group $\mathbb D_3$ generated by $\{g_1,g_2,g_3\}$. Note
that for Clairaut equation (\ref{Clairaut})  the above defined
three foliations   are $p=const$, $p+q=const$ and $q=const$
respectively.

\begin{figure}[th]
\begin{center}
\epsfig{file=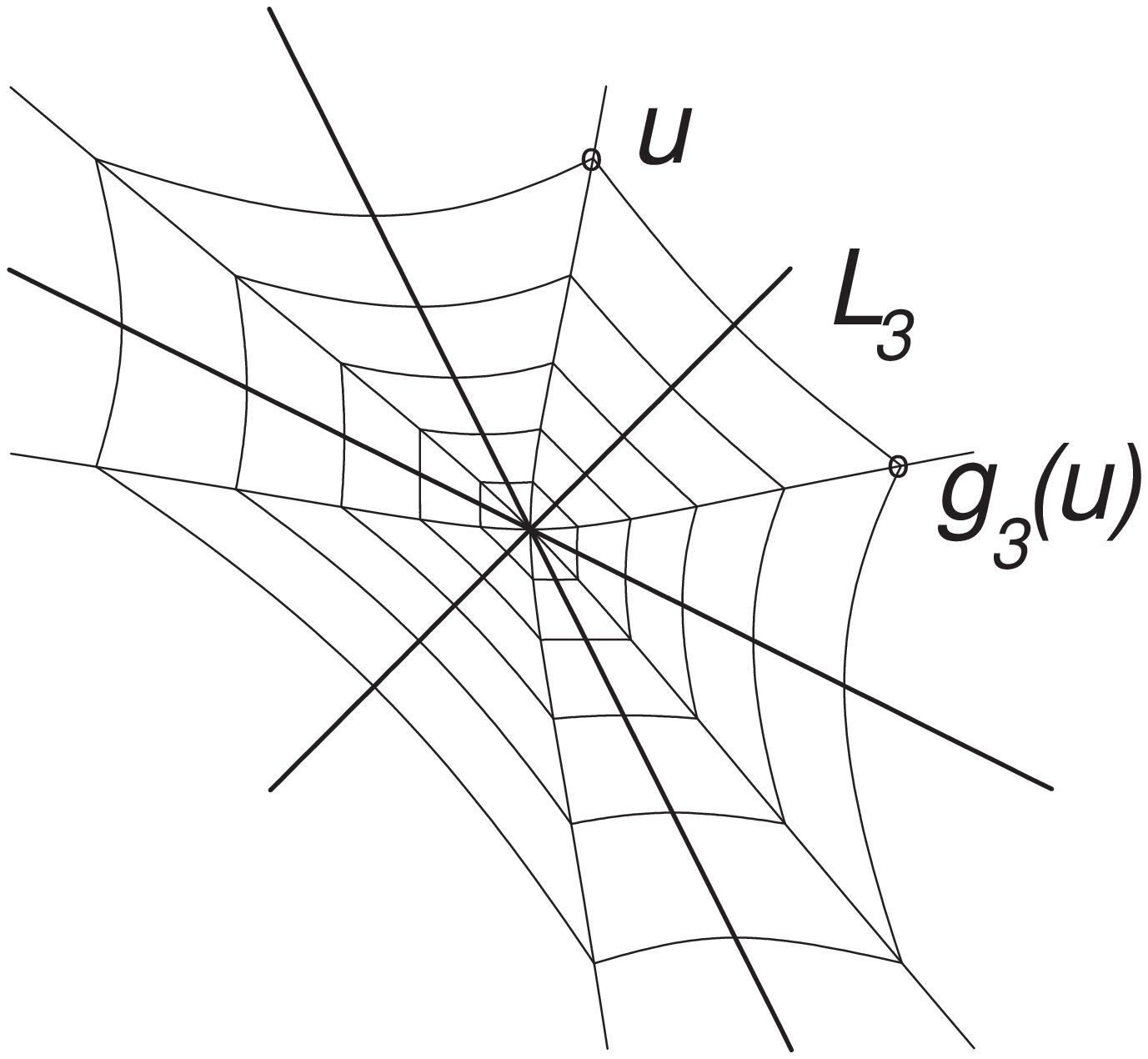,width=50mm}\ \ \ \ \ \ \ \ \ \ \ \ \ \ \ \ \
\ \epsfig{file=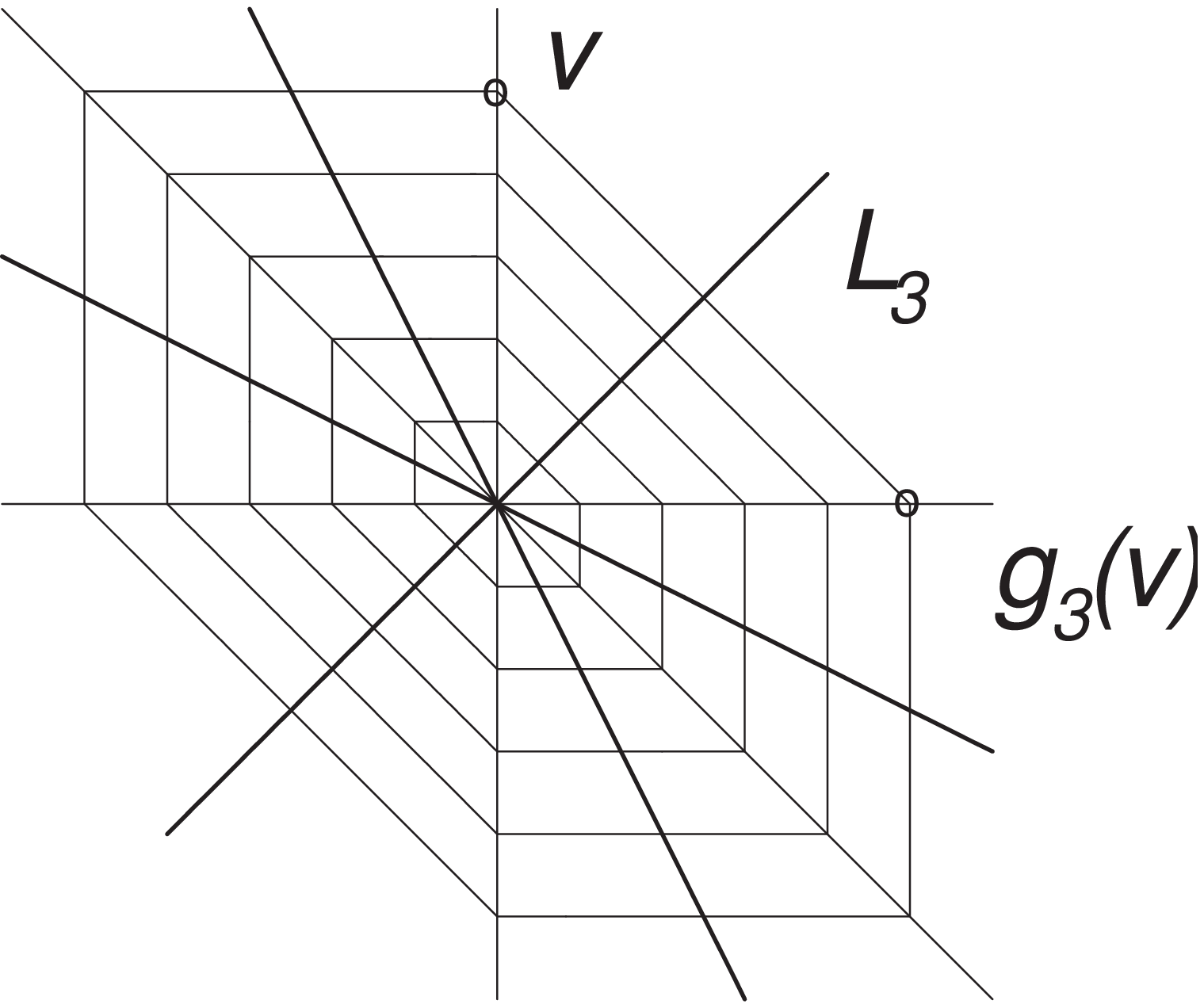,width=40mm}
  \caption{Brian\c{c}on's hexagons and their inverse images under the diffeomorphism $\psi$
with $u=\psi (v)$. } \label{Pic_web}
\end{center}
\end{figure}

Now we are ready to construct the diffeomorphism $\varphi$ that
transforms the given ODE to normal form (\ref{Clairaut}). Consider
a domain $U,0\subset \mathbb R^2,0$ such that the 3-web formed by
the foliations $\mathcal{F}_1$, $\mathcal{F}_2$ and
$\mathcal{F}_3$ is regular in $U$.  Let $\gamma _1$ be the
integral curve of $\tilde{\tau}$ that passes through the origin.
Pick up a point $u=(p_1,q_1)\in \gamma _1$ on this curve and draw
the {\it Brian\c{c}on} hexagon around $0$ through $u$. (Let us
recall the construction of the Brian\c{c}on hexagon: one draws
three curves $\gamma_i$, $i=1,2,3$ of the foliations
$\mathcal{F}_i$, $i=1,2,3$ through the origin, picks up a point on
one of this curves, say $\gamma_1$, and then goes around the
origin along the foliation curves, swapping the family whenever
one meets one of the $\gamma_i$. The web is hexagonal iff one gets
a closed hexagonal figure for any choice of the central "origin"
point and $u$. See Fig. \ref{Pic_web} on the left.) Let us choose
$u$ so that the following conditions hold:
\\1) $q_1>0$, \\
2) the Brian\c{c}on hexagon around $0$ through $u$ is contained in
$U$.

 Then there is a unique local homeomorphism
$\psi: \mathbb R^2,0 \to \mathbb R^2,0$ such that \\1) $\psi
(0,1)=u$,\\ 2) $\psi (1,0)=g_3(u)$, \\ 3) it maps the foliations
$p=const$, $p+q=const$ and $q=const$ to the foliations
$\mathcal{F}_1$, $\mathcal{F}_2$ and $\mathcal{F}_3$ respectively.
(See Fig. \ref{Pic_web}). In fact, the points $u$ and $g_3(u)$
lies on the same curve of the foliation $\mathcal{F}_2$ since the
involution $g_3$ is a symmetry of $\mathcal{F}_2$. Further, there
is a unique diffeomorphism, mapping the triangle
$(0,0),(0,1),(1,0)$ to the "triangle" $(0,0),u,g_3(u)$ formed by
the curves of the foliations $\mathcal{F}_1$, $\mathcal{F}_2$ and
$\mathcal{F}_3$ (see \cite{BB} p.15). This map is uniquely
extended to the whole hexagon. Moreover, the constructed
homeomorphism $\psi$ is equivariant with respect to the action of
$\mathbb D_3$ defined above. The map $\psi$ is smooth (analytic)
if the foliations $\mathcal{F}_1$,$\mathcal{F}_2$,$\mathcal{F}_3$
are smooth (analytic). Really, according to \cite{BB} p.155, there
exists a smooth (analytic) map, taking the foliations  $p=const$,
$p+q=const$ and $q=const$ to $\mathcal{F}_1$, $\mathcal{F}_2$ and
$\mathcal{F}_3$ respectively, and this map is uniquely defined by
specifying the inverse image of $u$. Thus this map should coincide
with the above homeomorphism $\psi$.

Now consider the  map:
$$
\pi \circ f \circ \psi: \mathbb R^2,0\to R^2,0,\  \ (P,Q)\mapsto
(A,B)=\pi(f(\psi(P,Q))).
$$
In coordinates it reads as
$$
A=\alpha (P,Q),\ B=\beta(P,Q).
$$
Observe that the above map is symmetric with respect to the action
of $\mathbb D_3$. Then by the results on smooth functions,
invariant with respect to finite group action, the functions
$\alpha$ and $\beta$ must depend only on the basic invariants of
the above $\mathbb D_3$ group action (see \cite{Gf}) and
\cite{Ss}):
$$
A=\tilde{\alpha} (X,Y),\ B=\tilde{\beta}(X,Y),
$$
where
\begin{equation}\label{Vieta}
X=-P^2-PQ-Q^2,\ \  Y=PQ(P+Q).
\end{equation}
 We claim that $\tilde{\alpha}
,\tilde{\beta}$ are the components of the searched for
diffeomorphism $\varphi$. To prove that consider the following
commutative diagram:
\begin{equation}\label{diagram}
\begin{array}{ccc}
  \mathbb R^2,0 & \stackrel{\psi}{\longrightarrow} &\mathbb R^2,0\\
  \downarrow V &   &\downarrow \pi \circ f\\
   \mathbb R^2,0 & \stackrel{\varphi}{\longrightarrow}&  \mathbb R^2,0\\
\end{array}%
\end{equation}
where $V$ is Vieta's map (\ref{Vieta}). Applying the same results
on symmetric functions to $\theta :=V\circ \psi ^{-1}$ we see that
the differentiable map $\theta$ is inverse to $\varphi $ inside
the "cusped" domain where our ODE has 3 distinct real solutions.
This completes the proof in the real analytic case. For the smooth
case we apply $\varphi$ and for the reduced equation we consider
the first integral of the  direction field $\tau$ that coincides
with $p$ on the part $M_3$ of $M$ that is projected to the domain
with three real roots. Further, it is easy to construct through
homotopy  the $\pi$-lowerable diffeomorphism $\varphi '$ of $M$
that is identity on $M_3$ and moves the integral curves of $\tau$
to that of (\ref{Clairaut}). The searched for diffeomorphism is
$\varphi '\circ\varphi$
\vspace{-3ex}\begin{flushright}$\Box$\end{flushright}

\subsection{The case of non-Legendrian criminant}

\begin{theorem}\label{cusped_th}
If the following conditions hold for implicit ODE (\ref{implicit}) at a point $m=(x_0,y_0,p_0)\in M$: \\
1) ODE (\ref{implicit}) has a hexagonal 3-web of solutions,\\
2) $p_0$ is the triple root of (\ref{implicit}) at $\pi(m)=(x_0,y_0)$,\\
3) the criminant $C$  is transverse to the contact plane field  in some punctured neighborhood of $m$ and ${\rm rank}((x,y,p)\mapsto (F,F_p))|_C=2$, \\
then it is locally equivalent to
\begin{equation}\label{cusped}
P^3+2PX+Y=0,
\end{equation}
within the domain, where (\ref{cusped}) has three real roots, if
$F$ is smooth, \\and in some neighborhood of $(0,0)$,  if $F$ is
real analytic.
\end{theorem}
{\it Proof:} We follow the proof scheme for Theorem
\ref{Clairaut_th}. Namely we consider the pull-back of the form
$dy-pdx$ to $\mathbb R^2,0$ by the Vieta map $\pi \circ f$, where
$f$ is defined by (\ref{pull_to_pq}), duplicate this pull-back
form by linear involutions (\ref{g2_g3}) and find a local
diffeomorphism of $\mathbb R^2,0$ matching "lifted" 3-web of our
equation and that  of (\ref{cusped}). The difference to the
previous case of Legendrian criminant is that now the web is
singular; each foliation $\mathcal{F}_i$ has a saddle singular
point at $f^{-1}(m)$. Therefore the classical results on hexagonal
3-web are not of much use to find the diffeomorphism "upstairs".
We construct it through a homotopy of the first integrals of the
corresponding foliations.

\vspace{1pt} \noindent$\bullet${\it Differential forms of the
foliations.}\\ By Lemma \ref{cubic} equation (\ref{implicit}) is
equivalent to (\ref{CUB}) with $A(0,0)=B(0,0)=B_x(0,0)=0$, $\left.
\frac{\partial(A,B)}{\partial(x,y)}\right |_0\ne 0$ and $$
A_x(0,0)\ne 0,\ \ B_y(0,0)\ne 0.
$$ Thus the
tangent plane $T_mM$ to the surface $M$ at $m=(0,0,0)$ is the
plane $y=0$ and the tangent line to the discriminant curve at
$(0,0)$ is $y=0$. Therefore one can choose $p,A$ as local
coordinates on the surface $M,m$ and $A,B$ as local coordinates on
$\mathbb R^2,0$. Thus $x=X(A,B),\ \ y=Y(A,B),$ where
$A=-p^2-q^2-pq$ and $p,q,r=-p-q$ are roots of (\ref{CUB}). As
easily follows from Theorem \ref{quadratic non_legendrian_web},
the kernel of the pull-back form $\pi^*(dy-pdx)$ is tangent to the
curve $D$ defined by (\ref{preimage_discriminant}). That means
that the kernel of the form
$$
f^*(\pi^*(dy-pdx))=(2p+q)(-Y_A+pX_A+qY_B-pqX_B)dp+(2q+p)(-Y_A+p(X_A+Y_B)-p^2X_b)dq
$$
is tangent to the line $L_1:=\{(p,q):\ p+2q=0 \}.$ Writing the
above form as $\omega_1:=(2p+q)P(p,q)dp+(2q+p)Q(p,q)dq$ with
suitable $P,Q$ and passing to the coordinates $q,\ s=2q+p$ one
gets
$$
\omega_1 =(2s-3q)Pds+(sQ-2(2s-3q)P)dq
$$
hence the tangency condition implies $P|_{s=0}=0$. By the Hadamard
lemma $P=s\tilde{P}$ hence one obtains
$$
\omega_1 =s\{(2s-3q)\tilde{P}ds+(Q-2(2s-3q)\tilde{P})dq\}.
$$
Now consider the expression for
$P=-Y_A+(s-2q)X_A+qY_B-(s-2q)qX_B$. One has $Y_A=0,$ $X_A\ne 0$
from $A_xB_y\ne 0$, $B_x=0$. Since $A$ is quadratic and $B$ is
cubic in $p,q$, the term $Y_A$ does not have linear terms in
$p,q$. This implies $\tilde{P}(0,0)\ne 0$. Using again the
condition on $L_1$ one obtains
$$
Q-2(2s-3q)\tilde{P}=s\tilde{Q}.
$$

Let as normalize the forms vanishing each on its own family of
solutions to satisfy $\sigma_1+\sigma_2+\sigma_3=0$:
\begin{equation}\label{sigmas}
\sigma_1=(q-r)(dy-pdx), \ \ \sigma_2=(r-p)(dy-qdx), \ \
\sigma_3=(p-q)(dy-rdx).
\end{equation}
As shown above the pull-back of $\sigma_1$ is
\begin{equation}\label{sigma_1}
\tilde{\sigma_1}=(2q+p)^2\{(2p+q)\tilde{P}dp+(2(2p+q)\tilde{P}+(2q+p)\tilde{Q})dq\}
\end{equation}

\vspace{1pt} \noindent$\bullet${\it Connection form.}\\ Following
\cite{Be} consider the area form
$$
\Omega:=\tilde{\sigma_1}\wedge
\tilde{\sigma_2}=\tilde{\sigma_2}\wedge
\tilde{\sigma_3}=\tilde{\sigma_3}\wedge
\tilde{\sigma_1}=(Y_AX_B-Y_BX_A)(p-q)^2(2p+q)^2(2q+p)^2dp\wedge dq
$$
and the connection form
$$
\gamma:
=h_2\tilde{\sigma_1}-h_1\tilde{\sigma_2}=h_3\tilde{\sigma_2}-h_2\tilde{\sigma_3}=h_1\tilde{\sigma_3}-h_3\tilde{\sigma_1},
$$
where $h_i$ are defined by $$d\tilde{\sigma_i}=h_i\Omega.$$ Using
(\ref{sigma_1}) on obtains by direct calculation that
$$
h_1=\frac{R_1}{(2p+q)^2(p-q)^2},
$$
where $R_1$ is a smooth function of $p,q$. Applying the cyclic
permutation $p\rightarrow q,\ q\rightarrow r,\ r\rightarrow p$,
one gets: $\tilde{\sigma_2}=(2p+q)^2\bar{\sigma_2}$,
$h_2=\frac{R_2}{(2q+p)^2(p-q)^2},$ where $R_2$ and
$\bar{\sigma_2}$ are smooth.  Therefore
$$
\gamma=\frac{ \bar{\gamma}}{(p-q)^2 },
$$
with a smooth form $\bar{\gamma}$. Observe that the connection
form $\gamma$ is symmetric with respect to the linear
transformation group $\mathbb D_3$, generated by $g_1,g_2,g_3$.
Therefore $\gamma $ is smooth, i.e. $(p-q)^2$ divides
$\bar{\gamma}$.

\vspace{1pt}

\noindent$\bullet ${\it Existence of first integrals.}\\ As the
web is hexagonal the connection form $\gamma $ is closed.
Therefore there exists a unique $\mathbb D_3$-symmetric function
$\mu$, satisfying
$$
d\mu=\mu \gamma , \ \ \mu(0,0)=1.
$$
Further, the forms $\mu \tilde{\sigma_i}$ are also closed, thus
defying functions $u_i$ by
$$
d(u_i)=\mu \tilde{\sigma_i}, \ \ u_i(0,0)=0,
$$
satisfying the following equation, which is equivalent to the
hexagonality of the web:
\begin{equation}\label{hex_identity}
u_1+u_2+u_3\equiv 0.
\end{equation}
Observe that the function $u_1$ is skew-symmetric with respect to
$g_1$:
\begin{equation}\label{skew_u1}
g_1^*(u_1)=-u_1.
\end{equation}
This follows from (\ref{sigmas}) and from the invariance of $\mu$.
Applying Hadamard's trick one estimates $u_1$ as follows:
$$
u_1(p,q)=\int _0^1 \frac{d}{dt}u_1(tp,tq)dt=\int _0^1 \left
(p\frac{\partial}{\partial p}u_1(tp,tq)+q\frac{\partial}{\partial
q}u_1(tp,tq)\right )dt=
$$
$$=(2q+p)^2\int _0^1
 \mu(tp,tq) \left (
 p(2p+q)\tilde{P}(tp,tq)+q(2(2p+q)\tilde{P}(tp,tq)+(2q+p)\tilde{Q}(tp,tq))
\right )t^3 dt.$$ Collecting similar terms, using
$\tilde{P}(0,0)\ne 0$ and integrating one has
\begin{equation}\label{u_1}
u_1(p,q)=(2q+p)^3((2p+q)\hat{P}(p,q)+q\hat{Q}(p,q)), \ \ {\rm
 where}
 \ \ \hat{P} (0,0)\ne 0.
\end{equation}

\vspace{1pt}

\noindent $\bullet${\it Properties of the first integrals.}\\ It
follows from Malgrange's Preparation Theorem that any smooth
function $F$ of $(p,q)$ can be represented in the form
\begin{equation}\label{Malgrange_form}
F(p,q)=F_0(A,B)+pF_1(A,B)+qF_2(A,B)+pqF_3(A,B)+q^2F_4(A,B)+pq^2F_5(A,B),
\end{equation}
with smooth functions $F_i$. In fact, the identities
$p^2=-pq-q^2-A$, $p^3=-pA-B$, $p^2q=-pq^2+B$, $q^3=-qA-B$ imply
$\langle p,q \rangle ^4 \subset \langle A,B \rangle$ and
$\mathcal{ E}(\mathbb R^2) / \langle A,B \rangle =\mathbb R
\{1,p,q,pq,q^2,pq^2 \}.$ (Here $\mathcal{ E}(\mathbb R^2)$ is the
local algebra of smooth map germs at $(0,0)$, $\langle p,q
\rangle$ its maximal ideal, generated by the coordinate functions
$p$ and $q$, $A: (p,q)\mapsto -p^2-pq-q^2$, $B: (p,q)\mapsto
p^2q+q^2p$, and $\mathbb R \{1,p,q,pq,q^2,pq^2 \}$ is the real
vector subspace of $\mathcal{ E}(\mathbb R^2)$, spanned by
$1,p,q,pq,q^2,pq^2$.) Moreover, inside the "cusped" domain with 3
real distinct solutions of our ODE the functions $F_i$ are
uniquely determined by $F$. For $F=u_1$ property (\ref{skew_u1})
implies $F_4=-F_3,$
 $F_2=2F_1-AF_5$, $F_0=-AF_3-BF_5/2$. Applying $-g_2^*$ and
 $-g_3^*$ to $u_1$ one gets the other two first integrals $u_2$
 and $u_3$, whose representations in form (\ref{Malgrange_form})
 are easily read from the representation of $u_1$.
Now identity (\ref{hex_identity}) implies $F_5=0$. Thus
\begin{equation}\label{symmetric_form_u1}
u_1(p,q)=(2q+p)(F_1(A,B)+pF_3(A,B)).
\end{equation}
Using (\ref{u_1}) one can write
$$
F_1(A,B)+pF_3(A,B)=(2q+p)^2G(p,q).
$$
Representing the function $G$ as
$$
G(p,q)=G_0(A,B)+pG_1(A,B)+qG_2(A,B)+pqG_3(A,B)+q^2G_4(A,B)+pq^2G_5(A,B)
$$
and substituting this representation into the above equation, one
obtains \begin{equation}\label{G} G(p,q)=
\frac{4}{3}AG_4(A,B)+pG_1(A,B)+pqG_4(A,B)+q^2G_4(A,B).
\end{equation}

 \vspace{1pt}

\noindent $\bullet${\it Equivariant homotopy.}\\ For equation
(\ref{cusped}) the first integral $u_1$ is $u_0(p,q):=p(2q+p)^3$.
Let us scale $u_1$ so that $\hat{P} (0,0) = \frac{1}{2}$ in
(\ref{u_1}). We claim that the family of functions
$$
u_t(p,q):=u_0(p,q)+t(u_1(p,q)-u_0(p,q))
$$
is equivariantly $\mathcal{R}$-trivial, i.e. for any $t\in [0,1]$
there is a diffeomorphism $\psi _t$, equivariant with respect to
above $\mathbb D_3$ group action, such that
$$
u_t\circ \psi _t=u_0.
$$
To prove this it is enough to find $\mathbb D_3$-equivariant
vector field $\xi (p,q,t) \partial_p+\eta (p,q,t) \partial _q$
satisfying the following homotopy equation:
$$
\xi (p,q,t) \frac{\partial u_t(p,q)}{\partial p}+\eta
(p,q,t)\frac{\partial u_t(p,q)}{\partial q}+\frac{\partial
u_t(p,q) }{\partial t}=0, \ \ \xi (0,0,t)=\eta (0,0,t)=0.
$$
A general form of a $\mathbb D_3$-equivariant vector field is
given by
\begin{equation}\label{equivariant}
\xi = p\alpha (A,B)+\left(\frac{A}{3}+pq +q^2\right)\beta (A,B), \
\ \eta = q\alpha (A,B)-\left(\frac{2}{3}A+q^2 \right)\beta (A,B),
\end{equation}
 (see \cite{Aw} or derive it from the
representations of $\xi, \eta$ in form (\ref{Malgrange_form})).
Observe that the difference $(u_1(p,q)-u_0(p,q))=\frac{\partial
u_t(p,q) }{\partial t}$ also has form (\ref{G}):
$$
u_1(p,q)-u_0(p,q)=\frac{4}{3}AL(A,B)+pK(A,B)+pqL(A,B)+q^2L(A,B)
$$
with $K(0,0)=0$ due to the chosen scaling of $u_1$.
 Solving the homotopy
equation yields the following expressions for $\alpha$ and $\beta
$:
$$
\alpha =
\frac{6K+(6K^2+2A^2(K_BL-KL_B)+9B(KL_A-K_AL)+10AL^2)t}{M},
$$
$$
\beta =\frac{-12L+(6A(KL_A-K_AL)+9B(KL_B-K_BL)+3KL)t}{M},
$$
where
$$
M=-24+[-48K-12AK_A-18B(K_B+2L_A)+8A^2L_B]t+
$$
$$
+[-24K^2-2A(6KK_A+25L^2)+3B(15K_AL-6KK_B-12KL_A)+
$$
$$
+2A^2(4KL_B-10LL_A-5K_BL)-30ABLL_B+(27B^2+4A^3)(K_AL_B-K_BL_A)]t^2.
$$
$M$ does not vanish at $(0,0)$ since $K(0,0)= 0$. The claim on
$\mathcal{R}$-triviality of the family $u_t(p,q)$ is proved.

\vspace{1pt}

\noindent $\bullet${\it Diffeomorphism.}\\ We have proved that the
diffeomorphism
$$
\psi := \psi _1
$$
maps the fibres of $u_0$, i.e. the curves $\{(p,q):
u_0(p,q)=const\}$ to those of $u_1$. Therefore, being equivariant,
$\varphi $ maps the foliations $\mathcal{F}_i$ of (\ref{cusped})
to that of our equation (\ref{CUB}). Now diagram (\ref{diagram})
defines again the desired diffeomorpfism $\varphi $. \hspace{230pt} $\Box$\\

\begin{figure}[th]
\begin{center}
\epsfig{file=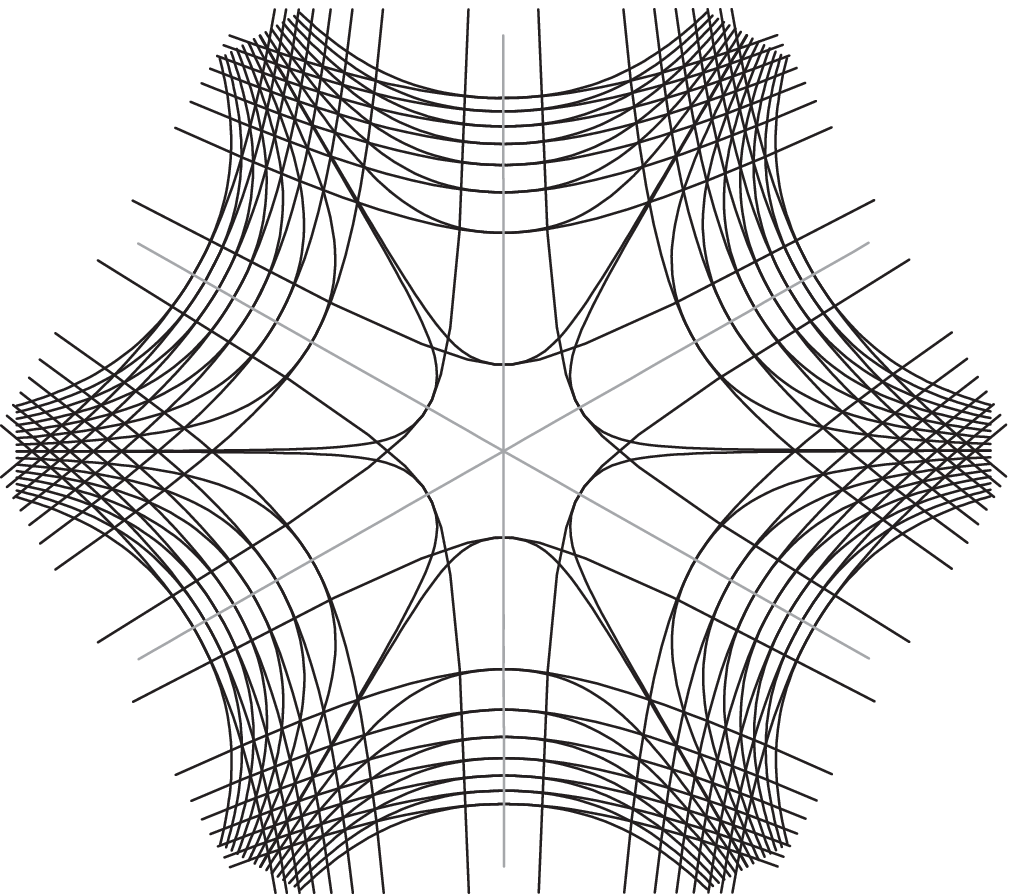,width=52mm}\ \ \ \ \ \ \
\epsfig{file=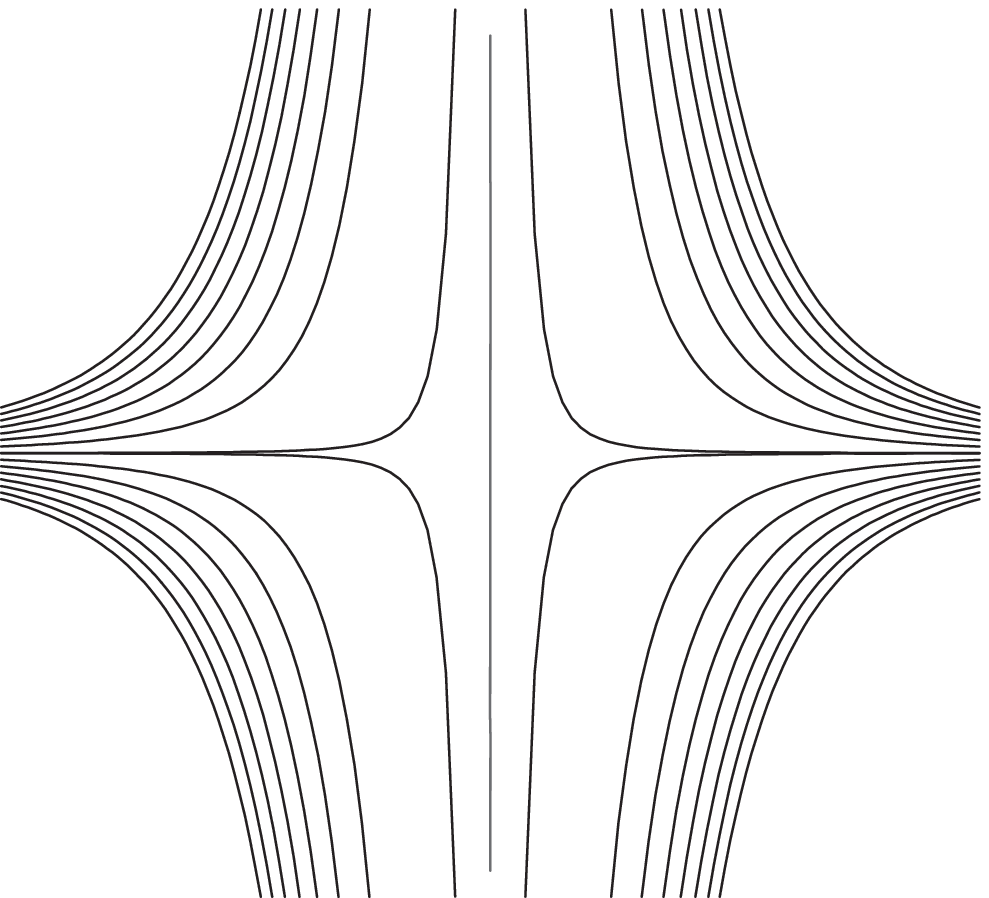,width=46mm} \ \ \ \ \ \ \
\epsfig{file=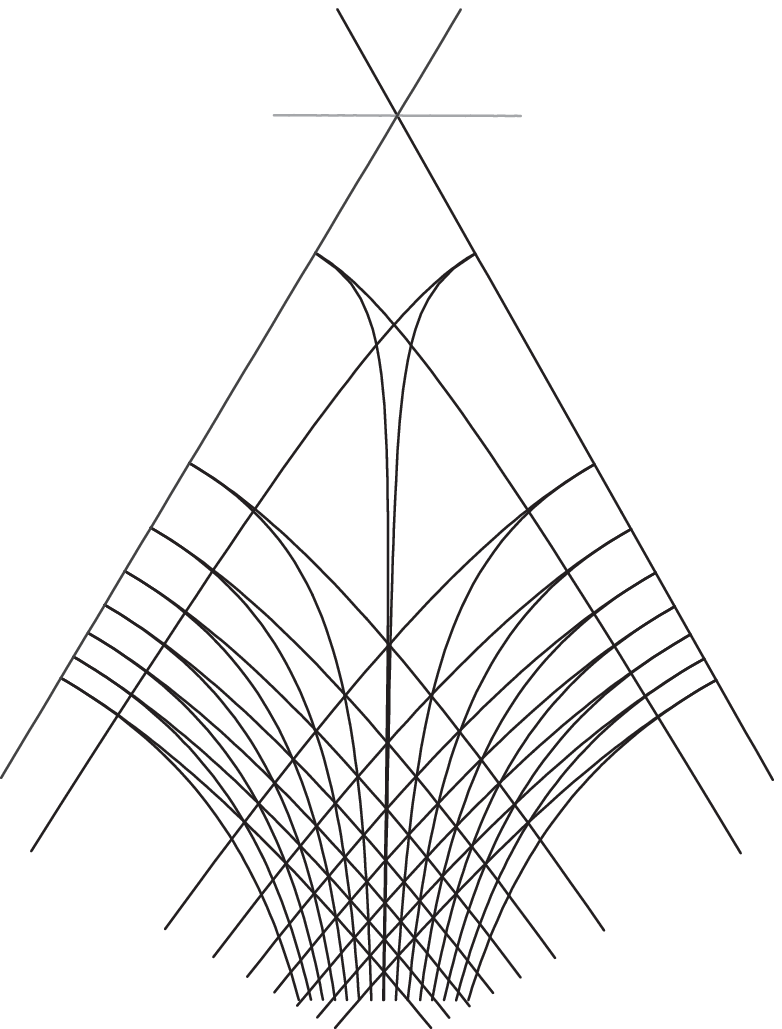,width=40mm}
  \caption{$\mathbb D_3$-symmetric hexagonal web of 3 foliations with saddle singularities (left),
   one of the foliations (center) and the fundamental domain of $\mathbb D_3$ (right).} \label{Pic_saddle}
\end{center}
\end{figure}

\noindent {\bf Remark.} The pictures of the $\mathbb
D_3$-symmetric hexagonal 3-web,  defined by the solutions of
(\ref{cusped}) and lifted to the plane $p+q+r=0$,  is presented in
Fig. \ref{Pic_saddle} on the left. It consists of 3 foliations,
one of them is shown in Fig. \ref{Pic_saddle} in the center. On
the right is the fundamental domain of $\mathbb D_3$-group
(compare with Fig. \ref{Pic_sol} on the right). The flower-like
form on the left suggests that the web is actually symmetric with
respect to the symmetry group $\mathbb D_6$ of regular hexagon. In
fact, it is the case since the fibers of the first integrals are
permuted by the following symmetry $g_4:\mathbb R^2 \to \mathbb
R^2,\ (p,q)\mapsto (-p,-q)$.

\section{Proof of the classification theorem}\label{proof}
Now we can prove Theorem \ref{classification_th}. If a point
$m=(x_0,y_0,p_0)\in M$ is regular then our equation is locally
equivalent to (\ref{normal}v) by definition. If $p_0$ is a double
root then the regularity condition (\ref{regularity}) implies that
$m$ is a fold point of the projection $\pi$. Futher Corollary
\ref{no_isolated_tangency} implies that the criminant is either
Legendrian or transverse to the contact plane field in some
neighborhood of $m$. Thus by Theorems
\ref{quadratic_legendrian_web} and \ref{quadratic
non_legendrian_web} the equation is locally equivalent either  to
(\ref{normal}iii) or to (\ref{normal}iv). Finally if $p_0$ is a
triple root and the criminant is either Legendrian or transverse
to the contact plane field  in some punctured neighborhood of $m$
then by Theorems \ref{Clairaut_th} and \ref{cusped_th} the
equation is locally equivalent either  to (\ref{normal}i) or to
(\ref{normal}ii). To complete the proof we show that Legendrian
and non-Legendrian parts of the criminant can not be glued
together at a cusp point. By Lemma \ref{cubic} our equation is
equivalent to (\ref{CUB}) with $(x_0,y_0,p_0)=(0,0,0)$,
$B_x(0,0)=0$. Suppose the criminant is transverse to the contact
plane field  for $p>0$ and Legemdrian for $p\le 0$. For any point
$m'\ne m$ with $p>0$ on the curve  $D$ defined by
(\ref{preimage_discriminant})   the  direction field $\tau$ is
tangent to $D$ by theorem \ref{quadratic non_legendrian_web}. This
condition reads as
$$
\left(\left. d(A+\frac{3}{4}p^2)\wedge (dy-pdx)\right)\right
|_{D}=0.
$$
In coordinates $p,x$ on $M$ it can be rewritten as follows
$$
\left(\left. (dA+\frac{3}{2}pdp)\wedge
((3p^2+A)dp+(B_x+(A_x+B_y)p+A_yp^2)dx)\right)\right |_{D}=0.
$$
Substituting $A=-\frac{3}{4}p^2$ and
$$
dA=A_xdx+A_ydy=A_xdx-\frac{A_y}{B_y+A_yp}((3p^2+A)dp+(B_x+A_xp)dx)
$$ into this equation one gets
$$
\left.
\frac{3}{2}p^2(A_xB_y+A_y(A_x+B_y)p+A_y^2p^2)-p(B_x+(A_x+B_y)p+A_yp^2)(B_y+A_yp)\right
|_{D}=0.
$$
Parameterizing the curve $D$ by $p$, expanding the above equation
by Tailor formula at $p=0$ and equating the coefficient by $p^2$
to $0$ one obtains
$$
\left. B_y(\frac{3}{2}A_x-(A_x+B_y))\right |_{x=y=0}=0,
$$
which implies
$$2B_y(0,0)-A_x(0,0)=0,$$
since $B_y(0,0)\ne 0$. (We have used the Taylor formula
$B_x=B_{xA}A+B_{xB}B+...=B_{xA}(-3/4p^2)+B_{xB}(-1/4p^3)+...$.)

On the other hand, for any point $m'\ne m$ with $p\le 0$ the
contact form vanishes on the criminant:
$$
\left. dy-pdx \right |_{C}=0.
$$
In coordinates $p,x$ on $M$ it can be rewritten as follows
$$
\left. (B_x+(A_x+B_y)p+A_yp^2)\right |_{C}=0.
$$
Now the Tailor expansion at $p=0$ for $C$ parameterized by p gives
$A_x+B_y=0$ ($B_x$ does not have linear in $p$ terms). Comparing
with the condition above on the non-Legendrian part one gets
$A_x(0,0)=0$ and therefore $B_y(0,0)=0$ which contradicts Lemma
\ref{cubic}.\vspace{-3ex}\begin{flushright}$\Box$\end{flushright}
\noindent {\bf Remark.} Unfortunately, the annoying stipulation in
Theorem \ref{classification_th} for the smooth case
(\ref{normal}i) can not be omitted to guarantee the existence of
the diffeomorphism $\varphi$ reducing ODE under consideration to
(\ref{normal}i) in some neighborhood of $m$ if one stays within
the framework of geometric Definition \ref{def_hex}. A necessary
condition for that is the existence of the first integral of
$\tau$ in the form
$f(p,x)^2g(p,q)^3=\tilde{\varphi}^*(p^2(x+3p^2/8)^3)$. Here
$\tilde{\varphi}$ is the lift to $M$ of the searched for
diffeomorphism and $p^2(x+3p^2/8)^3$ is the first integral of
$\tau$ for (\ref{normal}i). It is not hard to find a
counterexample which does not have such an integral in the form
$f(p,x)^2g(p,q)^3$ with not vanishing  $df,dg$ at $(0,0)$. This
drawback is repaired as follows. One replace definition
\ref{def_hex} with a less geometric one.
\begin{definition}\label{def_hex_PDE}
We say that implicit ODE (\ref{CUB}) has a hexagonal 3-web of
solutions if $A,B$ satisfy PDE (\ref{PDE}) and  the domain, where
(\ref{CUB}) has 3 real roots $p_1,p_2,p_3$ is not empty.
\end{definition}
The proof of Theorem \ref{cusped_th} is easily modified for the
case of one real root $p_1=p$ and two complex conjugated roots
$p_{2,3}=-p/2\pm iz$. The form $\sigma _1$ turns out to be pure
imaginary but the connection form $\gamma $ is real. All
analytical properties being the same, one finds the diffeomorpfism
$\varphi$ similarly through homotopy.

\section{Concluding remarks}\label{CR}

\noindent $\bullet$ {\bf Symmetries of the normal forms.} The
solutions of equation (\ref{normal}v) are the lines $dx=0$, $dy=0$
and $dx+dy=0$. The symmetry group of (\ref{normal}v)  is generated
by the following operators:
$$
X_1=\partial _x, \ \ X_2=\partial _y, \ \ X_3=x\partial
_x+y\partial _y.
$$
Thus the symmetry pseudogroup of a cubic implicit ODE with
hexagonal 3-web of solutions is at most 3-dimensional. In a
neighborhood of the projection of a regular point $m\in M$ it is
generated by the above three operators $X_i$ in suitable
coordinates. The coordinate change becomes singular on the
discriminant curve and not all symmetry operators "survive" at
$\pi (m)\in \Delta$. The symmetry pseudogroups of equations
(\ref{normal}iii) and (\ref{normal}iv) at a fold point are
generated by
$$
Y_1=x\partial_x+2y\partial_y,\ \ Y_2=\partial_x
$$
and
$$
Y_1=2x\partial_x+3y\partial_y,\ \ Y_2=\partial_y,
$$
respectively. This easily follows from Propositions
\ref{symmetry_quadratic_legendrian} and
\ref{symmetry_quadratic_non_legendrian}. Irreducible equations
(\ref{normal}i) and (\ref{normal}ii) have only one-dimensional
symmetry pseudogroup at $(0,0)$:
$$
Z=2x\partial_x+3y\partial_y.
$$

\noindent $\bullet$ {\bf Analytic properties.} All equations in
the given normal forms are integrable in elementary functions.

\noindent $\bullet$ {\bf Implicit cubic ODEs with singular
surfaces $M$.} Suppose our cubic ODE  factors out to 3 linear in
$p$ terms $p-f_i(x,y)$ such that 2 of 3 smooth surfaces
$M_i:=\{(x,y,p):p=f_i(x,y)\}$ intersect transversally along a
non-singular curve, the solutions of these 2 factors being
transverse to the curve projection into the plane. Then one can
bring these two factors to the forms $p=0$ and $p=2x$
respectively. The symmetry pseudogroup of the quadratic ODE
$p(p-2x)=0$ is $\tilde{y}=F(y), \ \ \tilde{x}=\sqrt
{F(y)-F(y-x^2)}.$ If our cubic equation has a hexagonal 3-web of
solutions then its third factor is generated by the vector field
$(\alpha(y-x^2)+\beta(y))\partial x+2x\beta (y)\partial _y$. As
the functions $\alpha ,\beta$ are arbitrary we can hope to "kill"
only one of them by the above mentioned symmetry. Thus a general
classification of all cubic ODEs will have functional moduli even
if one impose hexagonality condition. (Note that if the third
family of solutions in the example is transverse to the first two,
we have $\beta =-\alpha $ and one gets a finite classification
list.)

\noindent$\bullet$ {\bf Other examples.} The proof of Theorem
\ref{cusped_th} suggests the following procedure to generate cubic
ODEs with a hexagonal 3-web of solutions: start with a function
$F(p,q)$ written in form (\ref{Malgrange_form}) with
$A=p^2+pq+q^2,$ $B=pq(p+q)$, define $G(p,q):=F(q,-p-q)$,
$H(p,q):=F(-p-q,p)$ and solve the following  equations for $F_i$:
\begin{equation}\label{sym_hex}
g_3^*(F)=\pm F, \ \ F+G+H \equiv 0.
\end{equation}
This gives four of six coefficients $F_i$ as linear combinations
of the remaining two "free" functions of $A,B$. Then the fibers of
$F,G,H$ define a hexagonal 3-web, symmetric with respect to
$\mathbb D_3$-group action, generated by the involutions
$g_1,g_2,g_3$. The image of this web under Vieta map $(p,q)\mapsto
(p^2+pq+q^2,pq(p+q))$ is a hexagonal 3-web of solutions of some
implicit cubic ODE. For example, starting with $F=p-q$ one gets
the following equation:
$$
yp^3-\frac{2}{3}x^2p^2+xyp+\frac{1}{27}(2x^3-27y^2)=0.
$$
The solution 3-web of this equation is {\it dual} to that of
 (\ref{normal}i) (see \cite{Nw}). Its surface $M$ is not smooth at $(0,0,0)$.
  Note that one can also start with  $F$
 such that the "free" coefficients $F_i$ have poles at $0$. This
 approach linearizes the problem of finding local solutions
 of nonlinear PDE (\ref{PDE}). On the space of functions $F$ satisfying
 (\ref{sym_hex}) acts the pseudogroup of $\mathbb
D_3$-equivariant
 transformations with the tangent space generated by vector fields
 defined by
 (\ref{equivariant}). General classification of such functions with this equivalence
 group seems rather unpromising since the orbit codimension quickly becomes
 infinite.

\section{Acknowledgements}

The author thanks L.S.Challapa, M.A.S.Ruas, J.H.Rieger for useful
discussions. This research was partially supported by DAAD grant
415-br-probral/po-D/04/40407.

\end{document}